\documentclass[11pt,reqno]{amsart}

\pdfoutput=1
\numberwithin{equation}{section} \oddsidemargin=-.0cm
\usepackage{neuralnetwork}
\usepackage{amsmath}
    \usepackage{amssymb}
\usepackage{graphicx,scalerel}
\usepackage{bm} 

\usepackage{needspace}
\usepackage{newtxtext}
\usepackage{booktabs,caption} 

\usepackage[colorlinks=true, pdfstartview=FitV, linkcolor=blue,
            citecolor=blue, urlcolor=blue]{hyperref} 
\usepackage{subcaption}
\usepackage{sidecap}
\usepackage{xcolor}

\graphicspath{{./figures/},{../Figures/},{./Fig-QG/}}

\def\x{\boldsymbol{x}}

\definecolor{greenrb}{rgb}{0.2,0.6,0.2}

\definecolor{rred}{rgb}{0.7,0,0.1}
\definecolor{greenrb}{rgb}{0.2,0.6,0.2}
\newcommand{\mk}{\color{black}}
\newcommand{\mkr}{\color{rred}}

\def\bea{\begin{equation} \begin{aligned}}
\def\eea{\end{aligned} \end{equation}}
\def\beas{\begin{equation*} \begin{aligned}}
\def\eeas{\end{aligned} \end{equation*}}
\def\bes{\begin{equation*}}
\def\ees{\end{equation*}}
\def\be{\begin{equation}}
\def\ee{\end{equation}}

\newcommand{\Phib}{\boldsymbol{\Phi}}
\newcommand{\norm}[1]{\left\lVert#1\right\rVert}
\newcommand{\T}{\boldsymbol{\theta}}

\def\d{\, \mathrm{d}}

\newcommand{\Int } {\displaystyle \int}
\newcommand{\Frac}[2] {\frac{\textstyle #1} {\textstyle #2}}
\def\M{\, \mathcal{M}}
\def\R{\, \mathbb{R}}
\def\K{\, \mathcal{K}}
\def\L{\, \mathcal{L}}
\def\C{\, \mathcal{C}}
\def\D{\, \mathcal{D}}
\def\N{\, \mathcal{N}}
\def\G{\, \mathcal{G}}
\def\y{\mathbf{y}}
\def\z{\mathbf{z}}

\def\L{\mathcal{L}}

\def\adots{
  \mathinner{\mkern1mu\raise1pt\hbox{.}\mkern2mu\raise4pt\hbox{.}
  \mkern2mu\raise7pt\vbox{\kern7pt\hbox{.}}\mkern1mu}}
  
\newtheorem{thm}{Theorem}[section]
\newtheorem{lem}{Lemma}[section]
\newtheorem{defi}{Definition}[section]

\newtheorem{prop}[thm]{Proposition}

\newtheorem{rem}{Remark}[section]
\newtheorem{cor}{Corollary}[section]

\def\bt{\begin{thm}}
\def\et{\end{thm}}
\def\bl{\begin{lem}}
\def\el{\end{lem}}
\def\bd{\begin{defi}}
\def\ed{\end{defi}}
\def\bc{\begin{cor}}
\def\ec{\end{cor}}
\def\bp{\begin{proof}}
\def\ep{\end{proof}}
\def\br{\begin{rem}}
\def\er{\end{rem}}
\def\bi{\begin{itemize}}
\def\ei{\end{itemize}}

\title[Deep spectral computations]{Deep spectral computations in linear and nonlinear diffusion problems}

\author[E.~Simonnet]{Eric Simonnet}
\address[ES]{INPHYNI, UMR7010 CNRS-UNS, 1361, route des Lucioles, 06560 Valbonne}
\email{eric.simonnet@inphyni.cnrs.fr}

\author[M.~D.~Chekroun]{Micka\"el D.~Chekroun}
\address[MDC]{Department of Earth and Planetary Sciences, Weizmann Institute,
Rehovot 76100, Israel\\ and 
Department of Atmospheric and Oceanic Sciences, University of California, Los Angeles, CA 90095-1565, USA}
\email{mchekroun@atmos.ucla.edu}

\date{\today}

\begin{document}
\maketitle

\begin{abstract}
 We propose a flexible machine-learning framework for solving eigenvalue problems of diffusion operators in moderately large dimension. We improve on existing Neural Networks (NNs) eigensolvers by demonstrating our approach ability to compute (i) eigensolutions for non-self adjoint operators with small diffusion (ii) eigenpairs located deep within the spectrum (iii) computing several eigenmodes at once (iv) handling nonlinear eigenvalue problems. To do so, we adopt a variational approach consisting of minimizing a natural cost functional involving Rayleigh quotients, by means of simple adiabatic technics and multivalued feedforward neural parametrisation of the solutions. Compelling successes are reported for a 10-dimensional eigenvalue problem corresponding to a Kolmogorov operator associated with a mixing Stepanov flow. We moreover show that the approach allows for providing accurate eigensolutions for a 5-D Schr\"odinger operator having $32$ metastable states. In addition, we address the so-called Gelfand superlinear problem having exponential nonlinearities, in dimension $4$, and for nontrivial domains exhibiting cavities. In particular, we obtain NN-approximations of  high-energy solutions approaching  singular ones. We stress that each of these results are obtained using small-size neural networks in situations where classical methods are hopeless due to the curse of dimensionality. This work brings new perspectives for the study of Ruelle-Pollicot resonances, dimension reduction, nonlinear eigenvalue problems, and the study of metastability when the dynamics has no potential.
\end{abstract}

\section{Introduction}
In recent years, the idea of parameterizing the solutions to partial differential equation (PDEs) via a neural network (NN) has emerged as an influential approach to solve PDEs; e.g.~\cite{sirignano2018dgm,han2018solving,raissi2019physics,beck2021deep}. 
Unlike standard numerical methods which use meshes and thus are prone to the curse of dimensionality,  the partial derivatives of the NN approximation to the PDE solution, which, combined with the NN's natural ability  in representing high-dimensional functions, provides a powerful framework to overcome the curse \cite{han2018solving,beck2021deep}.  Still, challenges remain as each class of problems requires its own variational formulation that often takes into account prior knowledge about the problem's solutions through specific cost functionals and penalty terms \cite{ruthotto2020machine,nakamura2021adaptive}, along with its proper neural representation of the minimizers.

In spite of the recent compelling success of neural networks in representing high-dimensional PDE solutions with remarkable accuracy, the possibilities of NN-solvers for eigenvalues problems of differential operators  have been mainly explored in rather specific contexts \cite{han2020solving,yu2018deep,zhang2021solving}, although innovative ideas have emerged. For instance, Han et al.~\cite{han2020solving} propose to treat the eigenvalue problem of linear and semilinear second-order differential operators by reformulating it as a fixed point problem for the semigroup associated with the operator, exploiting the Feynman-Kac representation formula and forward-backward stochastic differential equations (FBSDEs). There, the eigenfunctions approximation is obtained via optimisation, through a neural-network ansatz, of a cost functional exploiting this representation. Operating in high dimensions, their algorithm allows for estimating the first eigenpair with an optional second eigenpair given some mild prior estimate of the eigenvalue.  

 In many works involving NN parametrizations, the 
diffusive coefficients are often large, yielding benchmark tests involving very smooth solutions.
Diminishing diffusive effects is notoriously difficult to handle especially in the context
of deep NNs. The problem is not so much related to the capacity of NNs to represent less regular solutions but rather to the difficulty of minimizing very stiff cost functionals.
We indeed show below that the expressivity of fully-connected NNs is often underestimated,
once the proper cost functional and minimization strategy are set. In particular, we show that 
nearly-singular solutions of diffusion problems---linear and nonlinear---can indeed be approximated with basic NNs.

In many applications such as described below, the computation of eigenmodes beyond the dominant ones is often very informative. Nevertheless, several intrinsic difficulties are tied to the eigenvalue problem of differential operators in high dimensions, calling for tempering the ambitions and focus to specific classes of problems.  
Among these difficulties one can mention the explosion of the number of eigenvalues with the dimension whose symptomatic behaviour is embodied by the asymptotic Weyl's formula describing the distribution $N(\lambda)$ of (large) eigenvalues {\mk for various diffusion operators. This is the case e.g.~of} the Dirichlet Laplacian in a bounded domain $\Omega$ of $\mathbb{R}^N$ \cite{weyl1911asymptotische,weyl1912asymptotische}, the Laplace-Beltrami operator \cite[Theorem 14.11]{zworski2012semiclassical} or {\mk Schr\"odinger operators \cite{arendt2009weyl,ivrii2016100,dyatlov2019mathematical}}. 
{\mk This combinatorial explosion is  the spectral signature of the curse of dimensionality. It may be furthermore amplified when a small parameter $\epsilon$ is present in front of the higher-derivatives; see \cite[Theorem 7.4]{dyatlov2019mathematical}.} 
For instance such a pathological behaviour is observed in the case for Schr\"odinger operators $-\epsilon \Delta + V(x)$, under certain assumptions on the potential,  for which the existence of $\sim \epsilon^{-n/2}$ resonances in specific bounded subsets of the complex plane is known to hold; see \cite{sjostrand1996trace,sjostrand2014weyl}. 

Aware of these difficulties, we focus in this article on three classes of diffusion problems, whose spectral investigations {\mk are} on a few but yet meaningful eigenpairs beyond the dominant ones.  Denoting by $\mathcal{M}$ either a smooth $N$-dimensional Riemannian manifold without boundary or the Euclidean space  $\mathbb{R}^N$ itself, our first focus is on the spectrum on Kolmogorov operators of the form  
\be\label{Eq_Kolmo_op}
\mathcal{K} \psi =\textrm{Tr}(Q (\x)D^2 \psi) +{\bf F (\x)} \cdot \nabla \psi, \; \x \in \M,
\ee
where $Q=(q^{ij})$ is a smooth mapping  from $\mathcal{M}$ taking values in the space of nonnegative symmetric $n \times n$ matrices, ${\bf F} = (F^i)$ is a (smooth) vector field on $\mathcal{M}$, and $\textrm{Tr}$ denotes the trace operator while ``$\cdot$'' denotes the inner product endowing $\M$.  

Our second focus is on the spectrum of Schr\"odinger operators of the form
\be\label{Eq_Schro_op} 
\mathcal{L} \psi=-\textrm{Tr}(Q (\x) D^2 \psi) + V(\x) \psi, \; \x \in \M,
\ee 
with $\M$ denoting the $n$-dimensional torus, and $V$, a potential function. 
For these problems, we are interested in situations where $Q$ is is taken constant and scales like $\epsilon >0$, with $\epsilon$ small.

Finally, our third focus is on the so-called nonlinear eigenvalue problems of the form
\bea\label{Eq_PDE}
-\Delta u & =\lambda f(u), \; \mbox{in} \; \Omega,\\
u&=0, \;\mbox{on} \; \partial \Omega,
\eea
where $\Omega$ is a bounded domain in $\mathbb{R}^N$ having a smooth boundary $\partial \Omega$, $\lambda$ is real, and $f$ is a superlinear positive function. Such problems have a long history and a strong mathematical basis \cite{rabinowitz1971some,amann1976fixed,lions1982existence,brezis1997blow} and arise in a wide range of disciplines, {\mk like in} gas combustion theory \cite{Gelfand63,bebernes2013mathematical,frank2015diffusion}, plasma physics \cite{chandrasekhar1957introduction,temam1975plasma}, or the theory of gravitational equilibrium of polytropic stars \cite{chandrasekhar1957introduction,fowler1931further,hopf1931emden}. Here the goal is to determine the multiplicity of solutions to \eqref{Eq_PDE} as $\lambda$ is varied, namely to {\mk compute} the bifurcation diagram.

Classical continuation or pseudo-arclength methods for computing the bifurcation diagram associated with Eq.~\eqref{Eq_PDE} may become challenging already in dimension 3.  This is for instance the case when the Laplacian of the solution grows in a superlinear, e.g. exponential, way. It thus calls for specific high-resolution local treatments and for handling the inversion of very large and ill-conditionned sparse matrices. An additional difficulty is due to the possible existence of an unbounded connected component of solutions bifurcating from infinity \cite{kielhofer2011bifurcation} associated with infinitely many turning points. Such situations are known to occur in dimension $3\leq N \leq 9$ for certain nonlinearities and domain geometries \cite{joseph1973quasilinear}. Obviously, the use of classical methods becomes hopeless as soon as $N\geq 4$.
Furthermore, many theoretical problems remain open in term of the dimension. For instance,  many theoretical problems are still open for $N\geq 3$, whether small bounded perturbations of $f$ can generate a discontinuity in the minimal branch through the appearance of a new fold-point;  see \cite[Theorem 3.1 \& Sec.~5]{chekroun2018topological}.

In this work, we propose a frontal approach for solving the eigenvalues problems described above, in a fully unsupervised way. The idea is to translate first the eigenproblems as minimization problems involving only the Rayleigh quotient together with ad-hoc normalisation constraints. We then parameterize the eigensolutions by simple feedforward NNs (FFNNs). In order to solve the optimisation problems, we adopt standard machine learning tools. We randomly sample points in the domain in order to estimate the cost functional gradients and perform a stochastic gradient descent until a statistical equilibrium is reached. 
  By doing so, the neural network learns the unknown function, bypassing the computational bottleneck inherent to grid-based methods.
 In the case of the Kolmogorov eigenvalue problem, we do not exploit simulation-based data for instance \cite{li2019computing}, or importance sampling technics  \cite{yan2022learning}. This means that in the context of e.g.~rare events calculation, no data is needed but only the knowledge of the governing equations. All the computational burden is on the minimization of ad-hoc cost functionals parametrized by  FFNNs. It is noteworthy that our approach is sufficiently general to handle other type of eigenvalue problems involving e.g.~higher-order derivatives, provided the eigensolutions have enough regularity.\\


 \section{Ruelle-Pollicott (RP) resonances in higher dimensions}

\subsection{Context}
The eigenvalues of the Kolmogorov operator \eqref{Eq_Kolmo_op} are also known as the Ruelle-Pollicott (RP) resonances \cite{Chekroun_al_RP2}.  These are encountered in many branches of physics (scattering resonances, statistical mechanics) and mathematics (zeta functions, dynamical systems); e.g.~\cite{ruelle1986locating,pollicott1986meromorphic,baladi2000positive,gaspard2005chaos,faure2011upper,lasota2013chaos,Cvitanovic2013,giulietti2013anosov,dyatlov2019mathematical}.

For deterministic systems when $Q=0$ in \eqref{Eq_Kolmo_op}, these resonances correspond to the eigenvalues of the transfer operator  \cite{schutte2001transfer} or its adjoint, the Koopman operator \cite{budivsic2012applied}.  They characterize fundamental properties of dynamical systems, such as power spectra,  mixing properties and decay of correlations \cite{pierrehumbert1994tracer,baladi2000positive,melbourne2007power,budivsic2012applied,lasota2013chaos,eisner2015operator}, coherent structures \cite{froyland2007detection,budivsic2012applied,froyland2013analytic,froyland2014almost},  metastability \cite{matkowsky1981eigenvalues,schutte2013metastability,pavliotis2014stochastic}, critical slow down \cite{tantet2015early,tantet2018crisis} or sensitivity to perturbations \cite{Chek_al14_RP,lucarini2016response,santos2020response}.  
Recently, RP resonances of stochastic systems \cite{Chekroun_al_RP2,RP_Hopf}  have shown their usefulness in the design of stochastic parameterizations to solve challenging closure  or data-driven model discovery problems issued from geophysical turbulence; see \cite{Chekroun2021c,santos2021reduced,KCB18}.

Different methods that have been developed over the last decades to compute RP resonances from finite-dimensional data-driven approximations of these infinite-dimensional operators, suffer from the curse of dimensionality. This is the case for instance of  methods rooted in the Ulam's approach. There, the underlying  transfer operator is approximated by Markov matrices giving an estimation of the transition probabilities from many short-term trajectories or a single long-term trajectory; see \cite{dellnitz1999approximation,schutte1999direct,schutte2001transfer,froyland2010coherent,Chek_al14_RP,tantet2015early,klus2018data,Chekroun_al_RP2,tantet2020ruelle}. Alternative approaches based on the infinitesimal generator while avoiding brute-force trajectory calculations, suffer from the same dimensionality restrictions as exploiting discretization and spectral collocation methods \cite{froyland2013estimating}.  

In parallel, the extended dynamic mode decomposition (EDMD) has been proposed as an alternative approach to approximate the spectral elements of the Koopman operator from multiple short bursts of simulation data \cite{williams2015data}.    
 The EDMD improves upon the classical dynamic mode decomposition (DMD) \cite{Rowley2009,Schmid2010} by the inclusion of a flexible choice of dictionary of observables to enrich the diversity of the spanning elements of the finite dimensional subspace from which the Koopman operator is approximated. Although the convergence of the EDMD has been established, applying the method in practice requires a careful choice of the observables to improve convergence towards Koopman's eigenfunctions with just a finite number of elementary bricks \cite{tu2014dynamic,williams2015data}.  
 
 This is especially difficult to achieve for high-dimensional and highly nonlinear systems.  In this case, the appropriate choice of observables  remains a challenge. Due to the presence of Koopman eigenfunctions with arbitrarily complex structures, it  
 may involve a large basis set to adequately represent them and their typical sharp gradients lying over the dynamics' separatices or unstable manifolds.  Such features call for dictionaries that are often manually curated, requiring problem-specific knowledge and painstaking tuning. 

More recently, iterative approximation algorithms which couples the EDMD approach with a trainable dictionary represented by an artificial neural network have been proposed to address this issue; see e.g.~\cite{li2017extended,yeung2019learning}. These machine-learning improvements of the EDMD enhances in essence the applicability of EDMD-based algorithms to approximate the spectral elements of the Koopman operator, from simulated data. The basic idea is to lift measurements to a higher-dimensional space where nonlinear problems tend to become more linear due to higher-dimensional embeddings (Cover's theorem \cite{mehrotra1997elements,lange2021fourier}).
Alternatively, autoencoder networks have been proposed to approximate Koopman eigenfunctions \cite{lusch2018deep}. The advantage is that of a low-dimensional latent space, which may promote interpretable solutions.

In spite of the great promises of these data-informed approaches and their recent deep learning directions, still open questions remain about how the choice of observables impacts the computation of the spectrum  \cite{brunton2021modern}. Even with the choice of ``good'' observables, the question of feeding the right regions of the phase space with the right amount of data constitute another practical barrier in applications. This is for instance the case of metastable systems perturbed by a small noise for which the proper sampling of rare events constitute an intrinsic challenge, especially in high dimension \cite{bouchet2019rare}.

 Instead, we present below a {\it simulation-free} approach, attacking directly the computation of eigenfunctions of the Kolmogorov operator without relying on data but rather exploiting its differential formulation whose coefficients depend on the governing equations. 

\subsection{Eigenmodes of $N$-dimensional Kolmogorov operators}
Thus, we consider Kolmogorov operators given in Eq.~\eqref{Eq_Kolmo_op} that we rewrite in coordinate form, 
\be\label{Eq_Kolmo_cf}
\mathcal{K}_\epsilon=\epsilon q^{ij} (\x) \partial_i \partial_j +F^i(\x) \partial_i, \; \partial_i=\partial/\partial x_i, 
\ee
 where the summation is taken over all repeated indices, and $\epsilon>0$ is a small parameter. 

To simplify the presentation, we restrict ourselves to the case of the $N$-dimensional torus  $\M=(-\pi,\pi)^N$, but our approach can be easily adapted to Kolmogorov operators on more general manifolds and with other boundary conditions.

We consider the Hilbert inner product $\left<a,b \right> = \int_{\M} a \bar  b \d \x$ with $L^2$-norm
$\norm{a} = \sqrt{\left<a,a\right>}$, and the Rayleigh quotient,
\begin{equation}\label{Rayq}
R_q = \frac{ \left< {\K}_\epsilon \phi,\phi \right>}{\norm{\phi}^2},~R_q \in \mathbb{C}.
\end{equation}
Solutions to the eigenvalue problem 
\be\label{Eq_eigpb}
\mathcal{K}_\epsilon \phi=\lambda \phi,
\ee
 satisfies $\norm{{\K}_\epsilon \phi - \lambda \phi} = 0$. A natural way is therefore 
to minimize this norm with respect to $\phi$ and $\lambda$. Although perfectly valid, we propose an alternative approach which is equivalent
and avoid having to handle in a separate fashion $\lambda$ and $\phi$ by noting that $\lambda = R_q(\phi)$, once $\phi$ is an eigenfunction.  In that respect, let us first remark that $\norm{{\K}_\epsilon \phi - \lambda \phi}^2 = \norm{{\K}_\epsilon \phi}^2 -|R_q|^2 \norm{\phi}^2$.
Naturally we also impose a norm constraint on the solution, here $\norm{\phi}=1$.
Let $\mu$ be a user-defined complex number. We then consider the following cost functional:
\bea\label{Eq_Cost}
{\C}_\mu(\phi) =  \gamma_{\rm cs} \overbrace{\left( \frac{\norm{{\K}_\epsilon  \phi}^2}{\norm{\phi}^2}
 - |R_q|^2 \right)}^{{\C}_s(\phi)} + \gamma_n ( \norm{\phi}^2-1)^2 \\
 + \gamma_{\rm bc} {\C}_{\rm bc} + \gamma_Q \left| R_q - \mu \right|^2.
\eea
The first term, ${\C}_s(\phi)$, is always positive by the Cauchy-Schwarz inequality and is zero if and only if $\phi$ is a solution of \eqref{Eq_eigpb}, 
 the second term normalizes the squared norm to be one, 
and the third term accounts for the boundary conditions. Here,   
\be
{\C}_{\rm bc} =
\sum_{\ell=1}^N \norm{\phi(\cdots,x_\ell,\cdots) - \phi(\cdots,x_\ell+2\pi,\cdots)}_{H^1(\mathcal{M})},
\ee
where $\norm{\cdot }_{H^1(\mathcal{M})}$ denotes the norm of the Sobolev space on the torus; see \cite[Chap.~9]{brezis_book}.
The last term $|R_q-\mu|^2$ constrains the eigensolution to stay in a particular user-defined region of the complex 
plane $\mathbb{C}$.  As mentioned above,  points in the domain are sampled randomly so that, in the course of
optimizing \eqref{Eq_Cost}, the solution does not depend upon a coordinate mesh.

The cost functional ${\C}_\mu(\phi)$ given in Eq.~\eqref{Eq_Cost} is minimized by means of  FFNNs (see Material \& Method) according to the following two-step procedure:
\begin{itemize}
\item[{\bf S1.}] {\bf Initial training of the NN}. It amounts to impose the NN-minimizer
to have (i) a large enough $L^2$-norm and (ii) a Rayleigh quotient close enough to the targeted
complex value $\mu$. This translates into having the coefficients $\gamma_n$ and $\gamma_Q$ to dominate the other penalty 
parameters.
\item[{\bf S2.}] {\bf Find the eigenpair $(\lambda,\phi)$ with $\lambda$ closer to $\mu$}. Once the NN has reached a statistical equilibrium, the NN-minimizer is relaxed by setting $\gamma_Q = 0$, while continuing the descent of the cost functional  ${\C}_\mu(\phi)$. This can be done abruptly like in this work or gradually. 
\end{itemize}
Step {\bf S1} must be thought as an optimal initialisation of the NN parameters for searching eigenfunctions having eigenvalues close to $\mu$.  It is therefore necessary to check whether the NN actually converges or not, given the state reached from Step {\bf S1}. It might happen for instance that the NN can drift away, sometimes to another eigenfunction, sometimes to some less relevant place of the landscape, e.g.~towards the trivial constant eigenmode with eigenvalue $\lambda=0$. It is therefore a good practice to reduce the training rate during Step  {\bf S2}.

\subsection{The multivalued deep learning of eigenstates}
 This approach extends naturally to the case of the simultaneous computation of multiple and distinct eigenpairs.
The generalization consists then simply to form a new cost functional obtained by summing up the cost functionals ${\C}_\mu(\phi)$ over a few targeted $\mu$'s. Steps {\bf S1} and {\bf S2} above are then followed, and at the end of the procedure, a single but multivalued FFNN is learned with  a  dimension output matching the number of targeted $\mu$'s. 

More precisely, we solve the following optimization problem
\be\label{Eq_multiopt}
\underset{\Phib(\cdot,\T)}\min \sum_{j=1}^N {C}_{\mu_j}([\Phib(\cdot,\T)]_j),
\ee
in which  $[\Phib(\cdot,\T)]_j$ denotes the $j$th component of the multivalued NN's output $\Phib(\cdot,\T)$; see Fig.~\ref{Fig_totoNN} in Material \& Method. 
Here, the real (resp.~imaginary) part of $[\Phib(\cdot,\T)]_j$ is aimed at approximating the 
real (resp.~imaginary) part of an eigenmode $\phi_j$ whose eigenvalue $\lambda_j$ is the closest to $\mu_j$. 
During the decent, the penalty parameters $\gamma_{Q_j}$ for which the $R_q$ have first converged,  are successively turned to zero.

Noticeable practical advantages are drawn from the usage of a single multivalued NN to compute several eigenmodes simultaneously. Indeed, not only this approach is simpler to implement than a counterpart that would consist of operating successively Steps {\bf S1}-{\bf S2} for distinct $\mu$'s, it leads to eigenmodes approximated in a much faster way.  The reason lies in a higher discriminant learning capability in this simultaneous, multiple-target setting, compared to the single-target setting. In fact, since in the single-target setting the convergence depends on the targeted eigenmode with convergence faster for certain modes than for others, in the multiple-target setting those that are found with fast convergence help constrain the FFNN to learn the others with less efforts.     
 
 This property can be interpreted as an intrinsic source of parallelism: it is likely faster to use a single NN with a multivalued output layer--one component per eigenmode to be approximated---rather than several independent scalar-valued NNs.
We do not know however what is a reasonable  upper bound of the size of the output layer, i.e.~how many eigenmodes can be computed   
simultaneously by this approach.

Finally,  we mention that such a multivalued deep learning of eigenmodes allows for avoiding fastidious grid search in the complex plane when e.g.~no spectral estimates are available by instead sampling randomly initial Rayleigh quotient values. 
The efficiency of such an approach is illustrated below on a 10-dimensional Kolmogorov eigenproblem.    

\subsection{RP resonances of $N$-dimensional stochastic mixing flows}
We address in this section the computation by our NN-solver of RP resonances associated with  a   
stochastic mixing flow on the 10-dimensional torus. The later is a stochastic perturbation of a 2-D Stepanov flow \cite{oxtoby1953stepanoff}, embedded within a 10-dimensional stochastic flow on the torus.  The 2-D embedded flow (variables $x_1$ and $x_6$) forces the other eight variables but not reciprocally; see  \eqref{10dflow} in Material \& Methods.  It is an instance of a {\it one-way coupled system}. Of course such a construction is somehow artificial but provides at the same time a dynamically challenging benchmark in terms of Kolmogorov spectrum. We explain why. 

 First, recall that  deterministic Stepanov flows are prototypes of flows that are topologically mixing on the torus \cite{oxtoby1953stepanoff}, that  exhibit already for 2-D flows a complicated temporal variability; see Fig.~\ref{Fig_Step_flows} below. {\mk In this case}, a non-trivial arrangement of the RP resonances is expected  in the complex plane \cite[Sec.~2.3]{Chekroun_al_RP2} associated with {\mk eigenmodes of singular structures exhibiting many sharp gradients distributed over thin confined regions, reminding those of a strange eigenmode \cite{pierrehumbert1994tracer}}.  Such complicated structures although smoothed out by noise, are expected to survive to a certain degree in  small-noise regimes. Already in 2-D, these structures are hard to approximate by Ulam's methods, requiring in particular a large amount of data to resolve the modes' fine structures  \cite{froyland2016optimal}.
       
\begin{figure}[htbp]
\centering
\includegraphics[width=.85\textwidth,height=.5\textwidth]{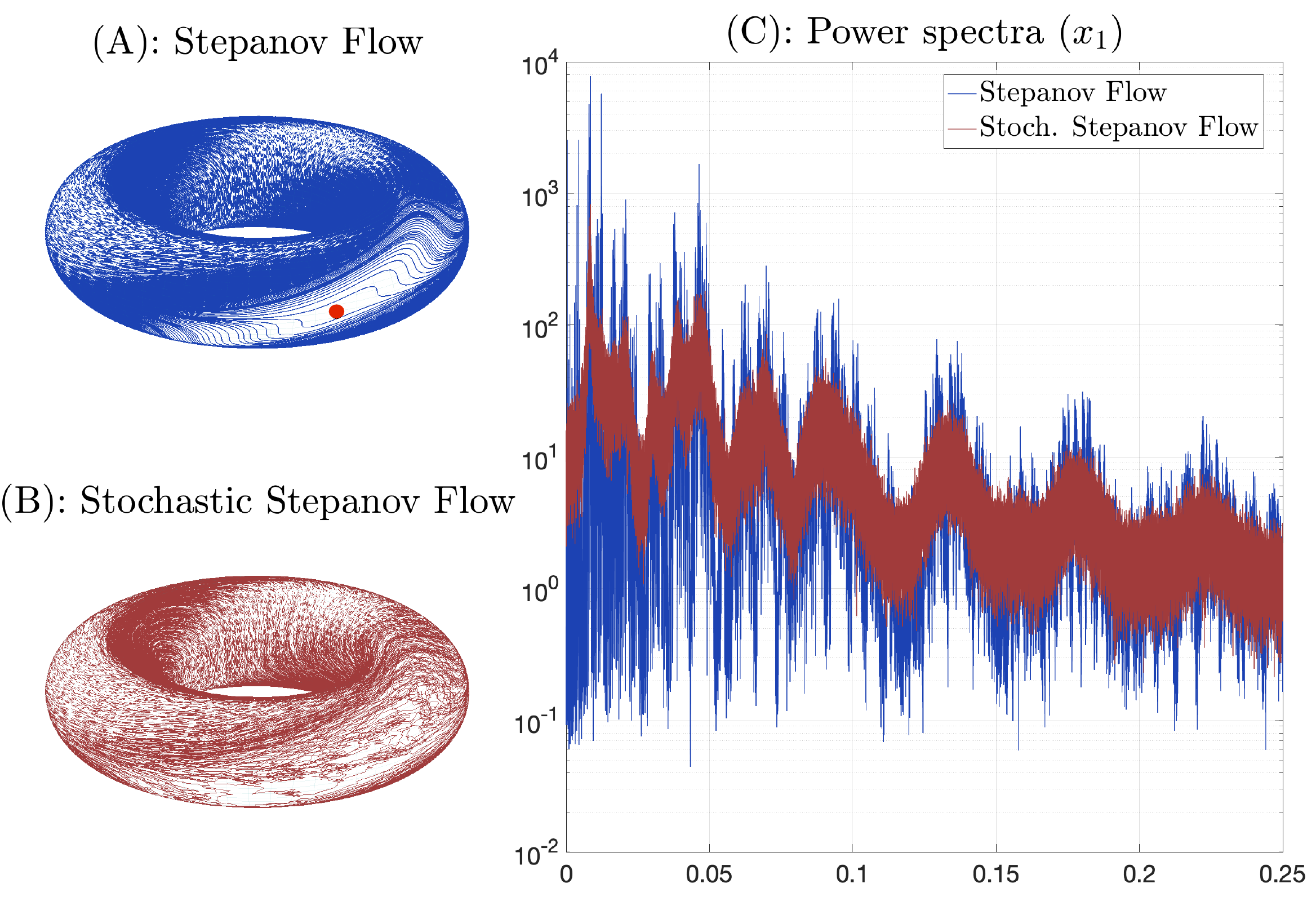}
\caption{{\bf Panel A}: A (segment of a) trajectory of the 2-D deterministic Stepanov flow solving \eqref{2dflow} on the torus. The unstable fixed point is shown by a red dot.
{\bf Panel B}: Same when \eqref{2dflow} is perturbed by an additive white noise $\epsilon (\dot{W_t}^2,\dot{W_t}^6)$, with $\epsilon=10^{-2}$. In both cases, $\alpha=\sqrt{20}$. {\bf Panel C}: Power spectra of these flows ($x_1$-variable); blue for deterministic flow, brown for the stochastic flow. Note that the noise smooths out the multiple bumps exhibited in the deterministic case.}
\label{Fig_Step_flows}
\end{figure}

These difficulties get severely amplified in dimensions higher than two. However, due to the one-way coupling in \eqref{10dflow}, a useful, low-dimensional characterization of certain  eigenmodes of the 10-D Kolmogorov operator associated with  \eqref{10dflow}, allows for testing our NN-solver's ability in resolving these issues. The proposition below summarises this point, whose proof is a simple exercise.

\begin{prop}  Let  $\y=(x_1,x_6)$ and $\z=(x_j)_{j\in I}$ with $I$ denoting the set of the first 10 positive integers to which $\{1,6\}$ is substracted. Let us write the Kolmogorov operator ${\K}_\epsilon$ associated with the 10-D stochastic system \eqref{10dflow} as
\be\label{Eq_LStepanov}
{\K}_\epsilon = \epsilon \Delta_{\y,\z} + {\bf F_1}(\y) \cdot \nabla_{\y}+ {\bf F_2}(\y,\z) \cdot \nabla_{\z},
\ee
where ${\bf F_2}(\y,\z)$  denotes the drift part in the RHS of \eqref{10dflow}  associated with the $\z$-variable and  where 
${\bf F_1}(\y)$ denotes the 2-D vector field  associated with the 2-D Stepanov system:
\be
\begin{array}{llll}\label{2dflow}
 \dot x_1  =  \alpha (1-\cos(x_1-x_6))+(1-\alpha)(1-\cos x_6)\\
  \dot x_6  =  \alpha(1-\cos(x_1-x_6)).
\end{array}
\ee
\end{prop}

Then the spectrum of ${\K}_\epsilon$ contains eigenfunctions of the type $(\phi^{2D}(\y),c)$, for any scalar $c$,  where $\phi^{2D}$ denotes any eigenfunction of the 2-D Kolmogorov operator 
\be\label{Eq_2DStepKolmo}
{\K}_\epsilon^{\y}=\epsilon \Delta_{\y} + {\bf F_1}(\y) \cdot \nabla_{\y}. 
\ee

\begin{figure*}[htbp]
\centering
\includegraphics[width=.9\textwidth,height=.55\textwidth]{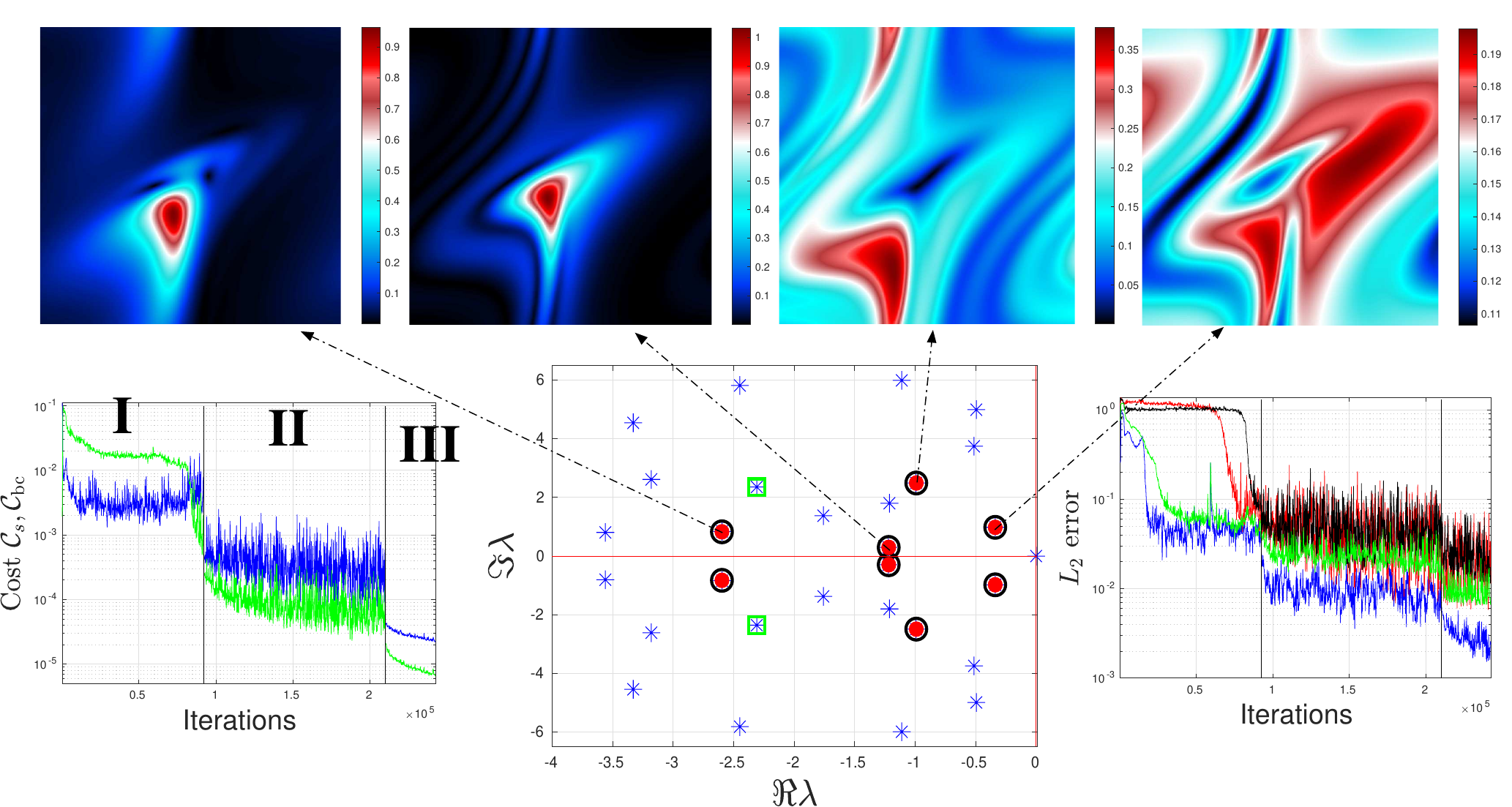}
\caption{Four eigenmodes of ${\K}_\epsilon$ given by \eqref{Eq_LStepanov} for $\epsilon = 0.1$ and associated with the 10D stochastic Stepanov flow \eqref{10dflow}. These are obtained by minimizing \eqref{Eq_multiopt} using simultaneous multivalued deep learning. Only the moduli are shown for a $L^2$-norm equals to one.
The central panel shows the benchmark eigenvalues (blue stars) and the ones obtained at the end of the learning process (red marks).
The lower-left panel shows the total cost behavior associated with $\gamma_{\rm cs}$ (blue curve) and $\gamma_{\rm bc}$ (green curve).
There are three stages. {\bf Stage I}:  $(\gamma_{cs},\gamma_n,\gamma_{\rm bc},\gamma_Q) = (15,2,10,10)$ imposes the 
NN to stay within some
prescribed region of the complex plane with learning rate $\alpha$, and batch size $n$ such that  $(\alpha,n) = (3 \cdot 10^{-3},1024)$. {\bf Stage II}: $(\gamma_{cs},\gamma_n,\gamma_{\rm bc},\gamma_Q) = (15,0.1,5,0)$ corresponds to a relaxation where the NN converges to the nearby eigensolutions with $(\alpha,n) = (10^{-3},1024)$.
{\bf Stage III} is a final ``precision run'' with
a smaller training rate $\alpha=10^{-4}$ and larger batches $n=4096$ (same penalisation coefficients).
The NN has 14 hidden layers with 25 neurons per layers using swish activation function, {\mkr resulting into 9,583 parameters}.
The right panel is the relative $L^2$-error compared with high resolution 2-D solutions evaluated
on each random batches. The mean relative $L^2$-errors over the last precision run 
are $0.3 \%$ (blue: leading mode shown to the upper-rightmost panel), $2.4 \%$ (black: low-frequency mode shown in the second upper-left), $2.6 \%$ (red: mode shown to the upper-leftmost panel) and $1 \%$ (green: mode shown in the second upper-right).}
\label{Fig_Stepanov}
\end{figure*}

The goal is thus to test whether our NN-solver is able to recognise the 2-D  embedded  patterns in the $\y$-variable exhibited by such eigenmodes of the 10-D Kolmogorov operator ${\K}_\epsilon$, in spite of the nonlinear coupling terms contained in ${\bf F_2}$.  In that respect, the 2-D patterns found by our NN-solver are benchmarked against the genuine 2-D eigenmodes obtained by solving the Kolmogorov eigenvalue problem  associated with \eqref{Eq_2DStepKolmo} using a standard method. Here, these 2-D modes are obtained over a $200\times 200$ finite-difference using a power iteration algorithm for $\epsilon=10^{-1}$. Finer resolutions are considered below for smaller $\epsilon$.

 The results shown in Fig.~\ref{Fig_Stepanov} demonstrate a striking success for the eigenmodes computation via minimization of \eqref{Eq_Cost} by using simple FFNNs. Not only the correct eigenmodes' patterns are found but also the two-dimensional feature of these modes are inferred, whether they are the dominant ones i.e.~close to the imaginary axis, or not.  In that respect, the underlying FFNN is able to identify on its own the essential variables governing the dynamics here $x_1$ and $x_6$, and the eigenmodes associated with these variables.  Such attributes are particularly relevant for dimensionality reduction, and will be discussed elsewhere. 
We focus next on another important challenge for applications, namely situations that are closer to the zero-noise limit.

\subsection{Approaching the zero-noise limit}
We illustrate here that our  framework  allows for the computation of eigenmodes close to the zero-noise limit, and located ``deep'' into the spectrum. 
Typically, the smaller $\epsilon$ is,  the harder the computation of such an eigenmode gets, even in low dimension, as mentioned above. The reason is that already for basic normal forms perturbed by noise (e.g.~pitchfork, Hopf) the limit is
singular and in many instances the Liouville eigenmodes ($\epsilon=0$) do not exist in a classical sense. They become Schwartz  distributions
and must be considered against smooth test functions/observables \cite{gaspard2001Hopf,gaspard2002trace,RP_Hopf}.
 The case of a 1-D pitchfork bifurcation is in that respect very informative. In this case, the eigenmodes are singular as involving the first derivatives of Dirac's distributions supported by the  unstable equilibria \cite{gaspard1995spectral}.

\begin{figure}[htbp]
\includegraphics[width=.85\textwidth,height=.45\textwidth]{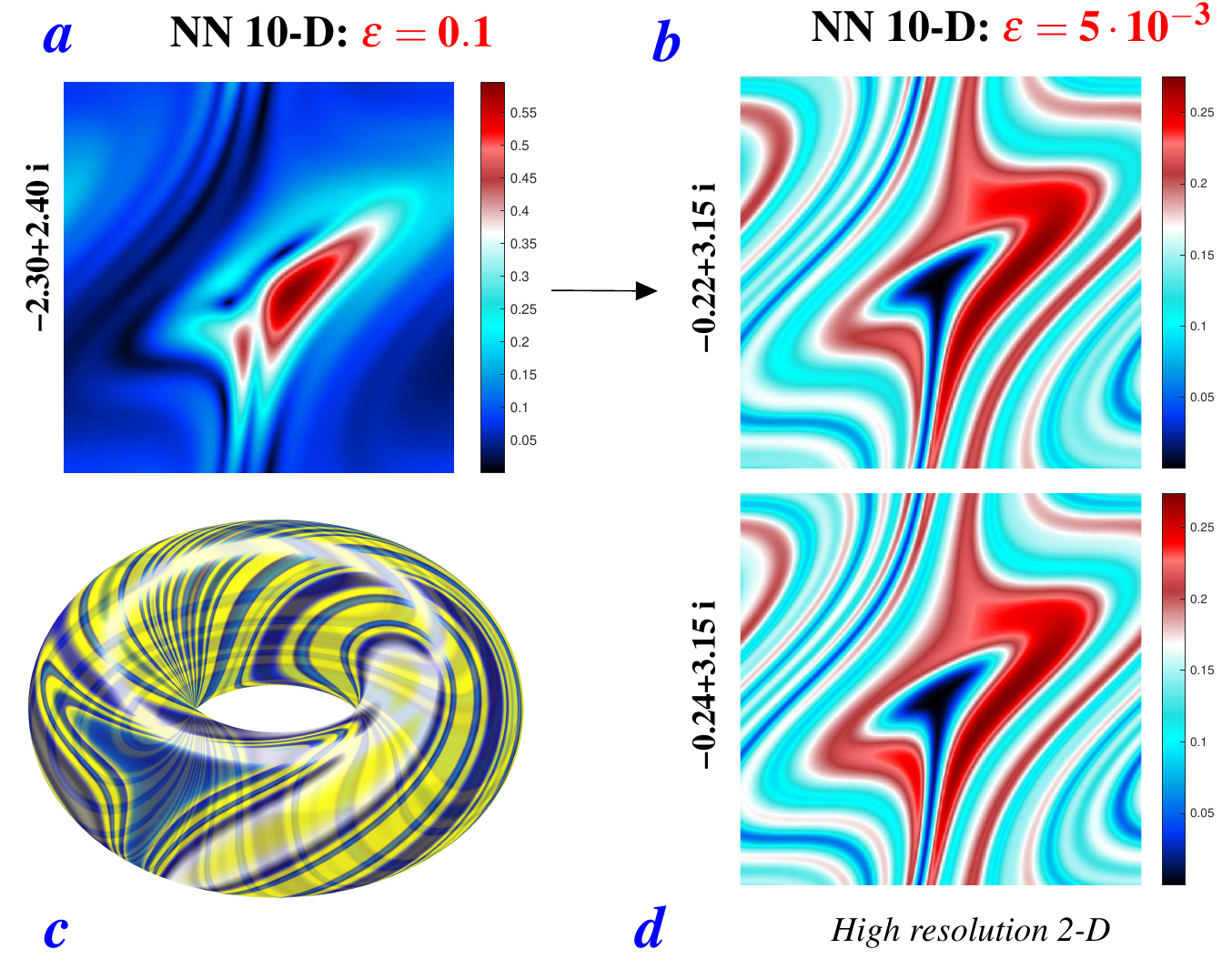}
\caption{ {\bf Adiabatic minimization of \eqref{Eq_Cost} using NN.} Starting from the eigenmode  shown in Panel (a)  for $\epsilon=0.1$ corresponding to the eigenpair marked in green in Fig.~\ref{Fig_Stepanov}, we slowly decreases
$\epsilon$ during the minimization of \eqref{Eq_Cost} on a logarithmic scale; typically in $50,000$ iterations. The eigenmode found this way by the NN for  $\epsilon = 5 \cdot 10^{-3}$ is shown in  Panel (b), for its modulus. 
The NN's architecture has $18$  layers, $30$ neurons,  with swish activation functions. During the precision run, the
batch size is $n=4096$ with a training rate $\alpha_r= 10^{-4}$.
The 2-D benchmark eigenmode has its modulus shown in Panel (d) as obtained by the MKL-Pardiso solver
using a $2000 \times 2000$-grid approximation  of ${\K}_\epsilon^{(\y)}$ given in \eqref{Eq_2DStepKolmo}. The corresponding eigenvalues are shown on the $y$-axis of Panels (a), (b) and (d). Panel (c) displays the mode shown in Panel (b) on the 2-D torus, using another color coding.}
\label{Fig_adia}
\end{figure}

For instance, by minimizing adiabatically \eqref{Eq_Cost} for the full 10-D Kolmogorov operator ${\K}_\epsilon$ with $\epsilon=5\times 10^{-3}$, our NN-solver is able to recover the 2-D embedded eigenmode's fine structures; cf.~Fig.~\ref{Fig_adia}-(d). The latter is obtained as eigenmode of the 2-D operator \eqref{Eq_2DStepKolmo} using a power iteration algorithm that exploits a high-resolution $2000\times 2000$ grid. Noteworthy is the much lower amount of parameters of the NN's architecture to achieve success here, namely about 17,000 parameters.

The thin and stretched structures encompassing a blue bulb-like pattern located around the unstable equilibrium in the center of Fig.~\ref{Fig_adia}-(b) are actually intimately related  to the topological mixing properties of the unperturbed flow. Although the deterministic Stepanov flow is ergodic, it has been indeed  numerically observed that a very long integration time is necessary for the dynamics to fill a small neighborhood of the unstable equilibrium (not shown). Over finite-time integrations, this phenomenon is accompanied by a dynamics' organization along ``strips'' of variable densities  (see Fig.~\ref{Fig_Step_flows}-(A)),  before reaching uniformity in the asymptotic limit. The mode shown in Fig.~\ref{Fig_adia}-(b) while located ``deep'' into the spectrum---corresponding to the green marks in Fig.~\ref{Fig_Stepanov}---is thus still very informative about the weak-noise limit.  As   $\epsilon$ is further decreased, the landscape exhibits sharper valleys leading eventually the NN to escape the  neighborhood of the targeted eigenmode.

\section{Eigenmodes of $N$-dimensional Schr\"odinger operators}
The next example we consider is inspired from \cite{han2020solving}. It consists of the following $N$-dimensional Schr\"odinger operator with periodic boundary conditions on the box ${\D}=[-2, 2]^N$, 
\be\label{Eq_Schro}
{\L}_{\epsilon} \psi = - \epsilon \Delta \psi + V(\x) \psi,
\ee
where the potential is given by 
\bes
V(\x) = \sum_{j=1}^N \left(-\frac{x_j^2}{2}+\frac{x_j^4}{4} + c_j x_j\right), 
\ees
with $c_j$ a scalar parameter. 

\begin{figure}[htbp]
\centering
\includegraphics[width=.85\textwidth,height=.4\textwidth]{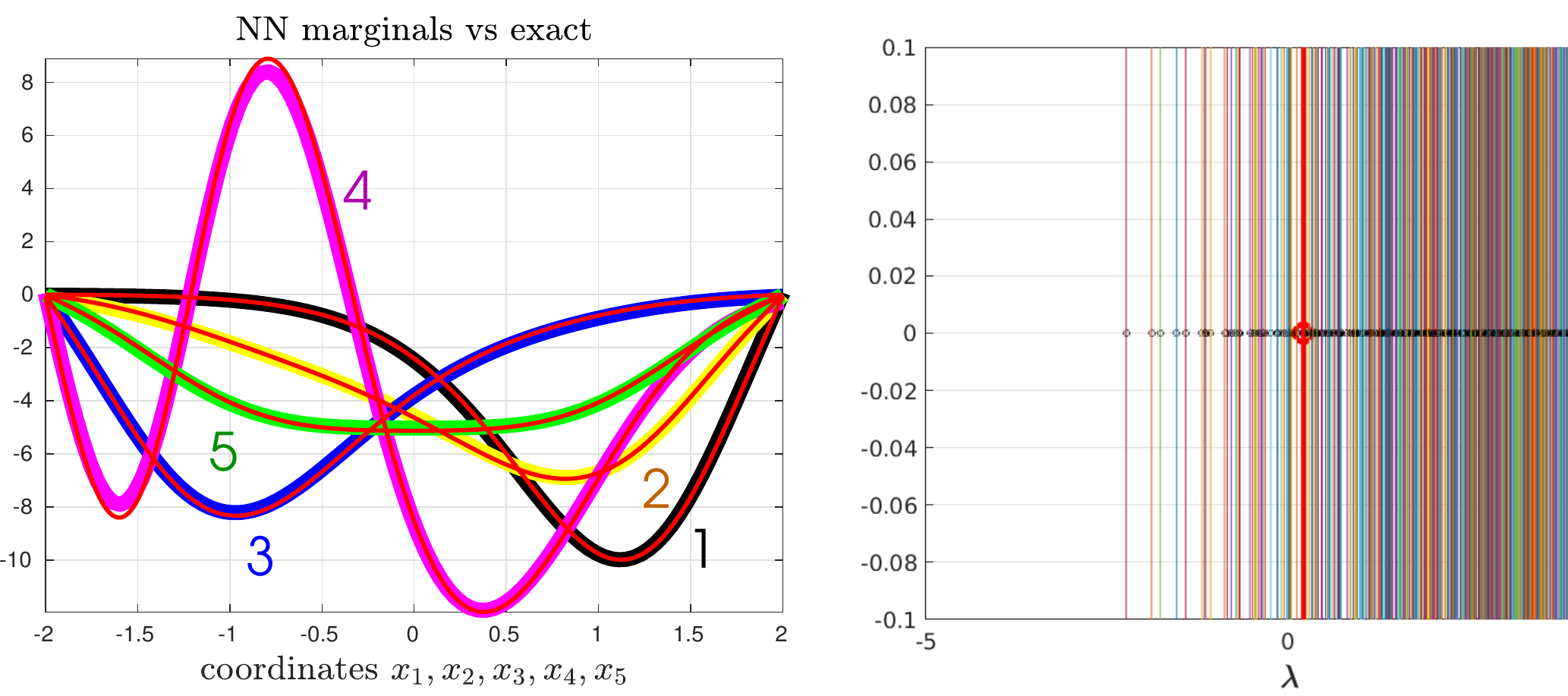}
\caption{A  5-D eigenmode (marginals shown) for the Schrodinger operator ${\L}_{\epsilon}$ given in \eqref{Eq_Schro} with $c=[\frac{1}{2}\; 1 \;\frac{3}{2} \; 2\; \frac{5}{2}]$ and $\epsilon = 0.2$,  ``deep'' into its eigenspectrum shown in the right panel, here the 45th-eigenvalue from the rightmost eigenvalue. 
The exact targeted eigenvalue is $\lambda= 0.203$ (red circle in right panel). Its approximation by our NN-solver is $\lambda_{\rm NN} =0.201$. The left panel shows a visual comparison of the NN marginals
$\bar{\phi}^{(k)}$ obtained from \eqref{Eq_marginal} with the $\psi_k$ (red curves) obtained from a high-precision finite-difference scheme. The underlying FFNN counts 12 hidden layers and 30 neurons per layer together with swish activation functions, resulting into 11,371 parameters. The batch size is $n=2048$ during the relaxation stage and the  ADAM's learning rate is $2\times 10^{-3}$.}
\label{Fig_Schro}
\end{figure}

The interest of this example is that the potential $V$ uncouples the variables, and is thus profitable for benchmark.   The full eigenvalue problem reduces indeed to solving $N$ disjoint eigenvalue problems for a 1-D Schr\"odinger operator, namely by solving $N$-times, $-\epsilon \psi_k''+(-x_j^2/2+x_j^4/4 + c_j x_k )=\lambda_k \psi_k$. The eigenvalues $\lambda$ of ${\L}_{\epsilon}$ are then obtained as sums of the $\lambda_k$, i.e.~$\lambda=\sum_{k=1}^N \lambda_k$, and the eigenmodes are given as product of the 1-D eigenmodes $\psi_k$, namely
\be\label{Eq_psi_full}
\psi(\x)=\prod_{k=1}^N \psi_k(x_k), \;\; \x=(x_1,\cdots,x_N) \in {\D}. 
\ee 
It is this tensorial property that makes interesting to submit to our NN-solver the $N$-dimensional  eigenvalue problem ${\L}_{\epsilon} \psi =\lambda \psi$. This way, one can test whether our framework allows for the NN to learn accurately the tensorial structure of the eigenmodes given by  \eqref{Eq_psi_full}. To do so, given an eigeinpair $(\lambda,\psi)$ obtained via \eqref{Eq_psi_full} and its NN-approximation $(\lambda_{\rm NN},\phi_{\rm NN})$,  it is thus sufficient to compare the $\psi_k$'s  in \eqref{Eq_psi_full} with the marginals $\bar{\phi}^{(k)}$ of $\phi$, given by
\be\label{Eq_marginal}
\bar{\phi}^{(k)}(x_k) =  \int_{\D} \phi_{\rm NN}(x_1,\cdots,x_n) \d{\x}^{(k)}, \; \d{\x}^{(k)} = \prod_{j=1,j \neq k}^N \d x_j.
\ee  
To test accuracy, the $\psi_k$'s are obtained by solving the corresponding 1-D Schr\"odinger eigenvalue problems,  using $1000$ grid points.  
 The results are shown in Fig.~\ref{Fig_Schro} for $N=5$, and for an eigenvalue that is located ``deep'' into the spectrum, namely  the $45$th eigenvalue of ${\L}_\epsilon$.

\section{Gelfand problem: High-dimensional bifurcations}
We conclude this article by considering the nonlinear eigenvalue problem, known  as the Gelfand problem, namely \eqref{Eq_PDE} with  $f(u)=e^u$ over a compact domain $\Omega$ in  $\mathbb{R}^N$. 
In the case of the unit ball $\Omega=B({\bf 0},1)$, due to the classical result of Gidas, Ni and Nirenberg \cite{gidas1979symmetry}, every solution to \eqref{Eq_PDE}  is radially symmetric and radially decreasing. The bifurcation diagram of the Gelfand problem---that provides the dependence on $\lambda$ of the solution set to \eqref{Eq_PDE}---is known to depend on the dimension, with in particular an infinite number of positive solutions for $\lambda=2(N-2)$, when $3\leq N \leq 9$; see \cite{joseph1973quasilinear}.   

In a first step, we benchmark the ability of our NN-solver to learn the radial symmetry of the solutions to \eqref{Eq_PDE} and the underlying bifurcation diagram with its first few turning points in dimension $N=3$ for the case $\Omega=B({\bf 0},1)$. The challenge is here, for the sake of generality, to do not rely on the radial symmetry which allows for transforming  \eqref{Eq_PDE} into the 1D-problem 
\be\label{Eq_1D}
u''=\frac{N-1}{r} u'+\lambda e^u, \;\; u(1)=u'(0)=0,
\ee
satisfied by the profile $u(r)=u(\norm{\x})$; see again \cite{joseph1973quasilinear}. 

Rather we aim at attacking the problem frontally in its original formulation \eqref{Eq_PDE}, to confront the ability of our NN-appraoch to handle the case of ``exotic'' geometries for which the transformation to \eqref{Eq_1D} does not apply. This is the case of domains exhibiting e.g.~cavities that break the symmetry for which much less is known theoretically with only partial results in special geometry like the annulus  \cite{nagasaki1994spectral}.

\subsection{The benchmark case: $\Omega=B({\bf 0},1)$ for $N=3$}
Traditionally, the bifurcation diagram for the Gelfand problem is shown in the $(\lambda,\norm{u}_\infty)$-plane. 
By the maximum principle \cite{gilbarg1977elliptic} every solution to this problem is positive. By the Gidas-Ni-Nirenberg symmetry result  \cite{gidas1979symmetry}, every solution $u_\lambda$ is radial and radially decreasing and thus its norm, $\norm{u}_\infty$, is attained at the center of the unit ball in the case $\Omega=B({\bf 0},1)$, i.e.~$\norm{u}_\infty=u({\bf 0})$.

It is well-known that there exists an extremal value $\lambda^\ast$ such that the nonlinear eigenvalue problem \eqref{Eq_PDE} has no solution, even in a weak sense for $\lambda >\lambda^\ast$; see \cite{brezis1996blow}. It is also well-known that the solution set $\{(\lambda,u_\lambda)\}$ forms an unbounded continuum in $[0,\lambda^\ast)\times C^2(\overline{\Omega}) $ that can be parameterized by a scalar $\tau$, with infinitely many turning points as $\lambda$ approaches the critical value $\lambda_s=2(N-2)$, $ 3 \leq N \leq 9$; see Material \& Methods. In this case, there exists for $\lambda=\lambda_s$ a singular solution $U_s(x)=-2 \log \norm{x}$; see \cite{joseph1973quasilinear}. More precisely, for $ 3 \leq N \leq 9$, the set of turning points, $T_k=(\lambda_{\tau_k},u_{\lambda_{\tau_k}})$, is infinite and $u_{\lambda_{\tau_k}}$ converges to $U_s$ in a weak sense as $k\rightarrow \infty$. As a consequence, the solution $u_{\lambda_{\tau_k}}$ takes large values and develops sharp gradients near the origin $\x={\bf 0}$ as $k\rightarrow \infty$, which makes extremely difficult the direct numerical computation of the bifurcation diagram by any standard continuation method for $N=3$ as one progresses across and above the turning points. Even if one uses the problem's radial symmetry and rely on the 1D-problem \eqref{Eq_1D} to compute the bifurcation diagram by a continuation method such difficulties survive as the solution's second derivative becomes exponentially large as $k\rightarrow \infty$.       

To handle such difficulties, Joseph and Lundgren \cite{joseph1973quasilinear} proposed an alternative two-step approach in which the boundary value problem \eqref{Eq_1D} is  treated by a shooting argument combined with an Emden's transformation facilitating a phase plane analysis to infer the bifurcation diagram; see Material and Methods.  This approach does not extend however to situations in which the domain's symmetry is broken such as considered below.

Thus, to address the bifurcation diagram computation for general situations within a variational approach suitable to an NN-treatment, we proceeds as follows. 
A first idea is to set a target value $A$ of the norm $\norm{u}_\infty$ and find the corresponding $(\lambda,u_\lambda)$ by minimizing the cost functional 
\be\label{Eq_initialC}
{\C}(u) = \rho_g \int_{\Omega} |\Delta u + \lambda {\rm e}^u|^2 \d \x  + \rho_{0} |u({\bf 0}) - A|^2 + {\rm B.Cs},
\ee
where $\rho_g $ and $\rho_{0}$ are positive free coefficients. 
Denoting by $\langle \cdot, \cdot\rangle$ the $L^2$-inner product, the eigenvalue to be found is then $\lambda = -\langle \Delta u_\lambda,\varphi \rangle / \langle {\rm e}^{u_\lambda},\varphi \rangle$ where $\varphi$ can be any reasonable test function (e.g.~$\varphi = 1$). The free parameter $A$ controls the energy level in the L$^\infty$-norm  of the solution $u_\lambda^A$  that the NN is aimed at approximating.  
This parameter plays a similar role than the constraint on the $L^2$-norm for the Stepanov eigenvalue problem, i.e.~the $\gamma_n$-penalty term in \eqref{Eq_Cost}. Also, by setting a few distinct $A$-values,  a multivalued approach may be adopted to approximate the corresponding ``eigenpairs'' $(\lambda,u_\lambda^A)$. Nevertheless, the obtention of good approximations of $u_\lambda^A$ by minimization of \eqref{Eq_initialC} is becoming more and more challenging as $A \rightarrow \infty.$

Indeed, the presence of a log-singular solution causes the second-order derivatives of any regular solution $u_\lambda^A$ to scale as $O({\rm e}^{A})$ as $\norm{u_\lambda^A}_{\infty}=A \to \infty$, near the ball's center. This phenomenon manifests into a  saturation of the NN ability to approximate correctly these derivatives, and {\it in fine} the $u_\lambda^A$ with high energy.  The problem's stiffness encountered in the direct computation of the bifurcation diagram via a continuation method for  \eqref{Eq_PDE} (or  \eqref{Eq_1D}) is here transposed into the minimization of an ill-conditioned problem \eqref{Eq_initialC}.        

Thus, the idea to revise \eqref{Eq_initialC}  by relying on an Emden-type transformation to ``blow-out'' the log-singularity.  This idea parallels Joseph and Lundgren's approach \cite{joseph1973quasilinear} albeit in a more general setting; see Material \& Methods. The change of coordinates we retain here is of the form   
\be\label{Eq_transfo}
{\bf X} = \frac{g(\norm{\x})}{\norm{\x}} \x,
\ee
where $g$ is defined by $g(r)=c \log \left(\frac{r}{\epsilon} + 1\right)$ with $\epsilon>0$ is small, and $c$ is chosen such that  $g(1)=1$, namely $c = 1/\log (1+\epsilon^{-1})$. Note that $g$ is invertible and given by $r = g^{-1}(R) = \epsilon \left( {\rm e}^{R/c} - 1 \right)$. This transformation allows us to map the unit ball onto itself and to alleviate the singular behavior at $R=\norm{{\bf X}}=0$ in the minimization of the revised cost functional; see \eqref{Eq_rev_cost} below.

To rewrite the cost functional \eqref{Eq_initialC} in this new coordinate system,  one needs to express the corresponding Laplacian which takes the form
\be\label{Eq_Deltanew}
\Delta_{\bf X}= \sum_{k=1}^N \Delta X_k \frac{\partial}{\partial {X_k}} + \sum_{i,j = 1}^N
\nabla X_i \cdot \nabla X_j \frac{\partial^2}{\partial {X_k}^2},
\ee
where $\Delta X_k$ and $\nabla X_i$ involve partial derivatives with respect to $\x$, but expressed in the new coordinates ${\bf X}$.
The coefficients in \eqref{Eq_Deltanew} are not radial and are moreover singular at ${\bf X}=0$. Exact
expressions of these coefficients are given in Material \& Methods.

\begin{figure}[htbp]
\vspace{-.2cm}
\includegraphics[width=.9\textwidth,height=.8\textwidth]{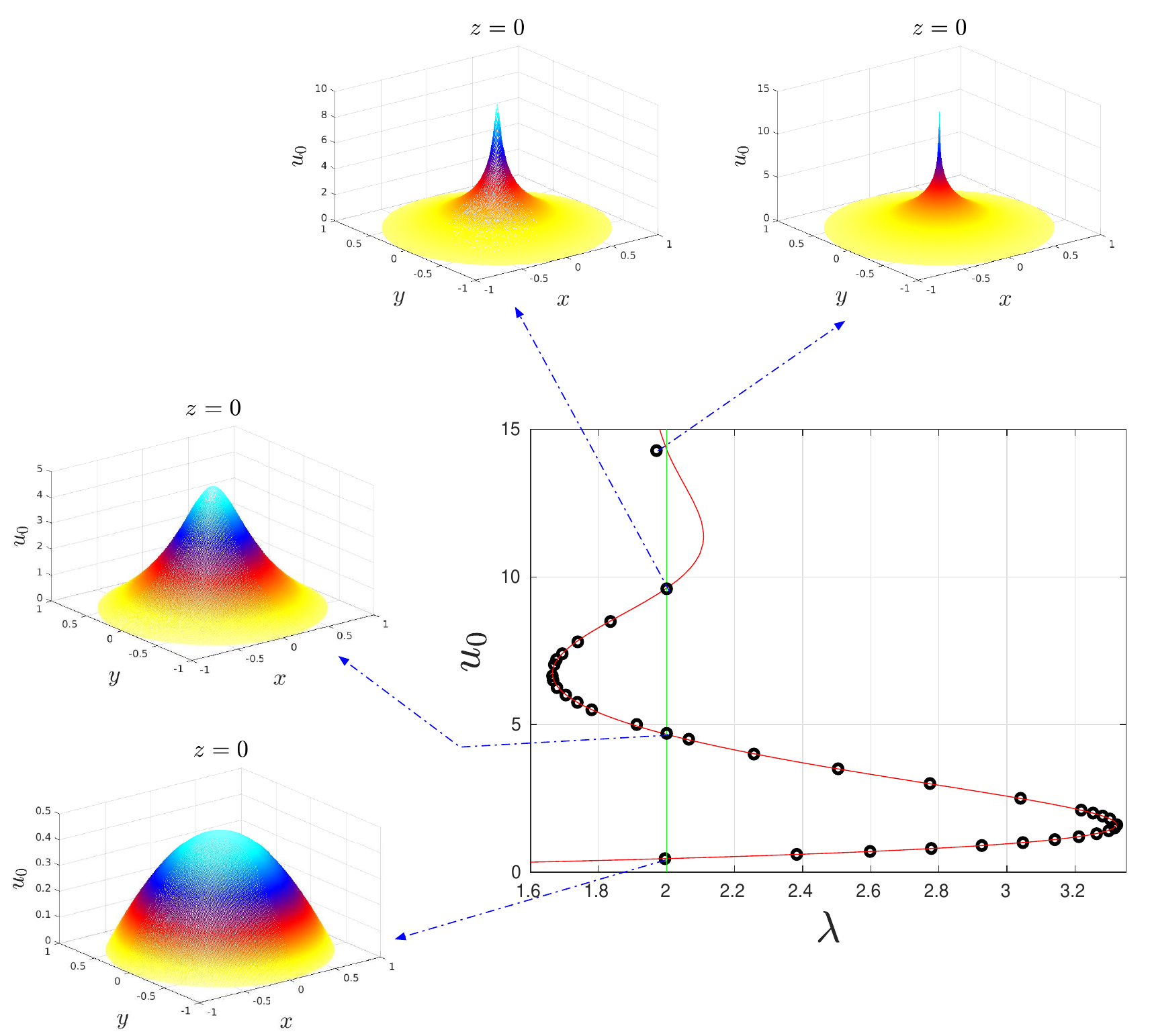}
\caption{Bifurcation diagram obtained in 3-D by minimization of \eqref{Eq_rev_cost} using a deep NN, following the steps (G$_1$) and (G$_2$) (black circles). 
Along this bifurcation diagram are shown cross-sections in the plane $z=0$ of four solutions obtained this way at the critical value $\lambda = 2$ ($N=3$), from low- to high-energy. The benchmark red curve is obtained by using a high resolution 10,000-grid pseudo-arclength code solving the 1-D (radial) problem \eqref{Eq_1D}. See Material \& Methods for the NN's configuration.}
\label{Fig_Gelfand}
\end{figure}

The Dirichlet boundary conditions on $S^2$ are handled using a simple lift idea. We thus write $u(X) = A \Theta(\norm{{\bf X}}) v({\bf X})$with $\Theta(R) = \cos (\pi R /2)$for instance. The function $v$ is the one parametrized by our NN. The usage of this lift makes  the problem unconstrained. Combined with the inverse of the transformation \eqref{Eq_transfo} it leads us finally to revise the minimization of \eqref{Eq_initialC} into the minimization of:
\bea\label{Eq_rev_cost}
&{\G}_A(v) =  \rho_g \Int_{\Omega} \left| \Delta_{\bf X} (\Theta v) +  \lambda {\rm e}^{A \Theta v} \right|^2  w^2(R) \d{\bf X}\\
&\hspace{3.3cm}+ \rho_{0} |v({\bf 0}) - 1|^2, \\
& \mbox{with } \lambda = -A \displaystyle \frac{\langle \Delta_{\bf X} (\Theta v), \varphi \rangle}{ \langle {\rm e}^{A \Theta v}, \varphi \rangle}, \; \varphi({\bf X})=w(\norm{\bf X}).
\eea
Here, the function $w$ is chosen to be $w(R) =g^{-1}(R)$. Our minimization of  \eqref{Eq_initialC} is then organized in two consecutive steps: 
\begin{itemize}
\item[(G$_1$)] Fix $A > 0$ and set $\rho_0 \gg \rho_g$ to enforce the NN-approximation ${\N}$ to satisfy ${\N}({\bf 0})\approx 1$, while ${\G}_A \ll 1$. 
\item[(G$_2$)] Set $\rho_0 = 0$: We relax the NN-approximation to solve only the nonlinear eigenvalue problem.
\end{itemize}
This approach allows us to compute a large portion of the bifurcation diagram with high precision.  The approach can be understood as a poor-man continuation approach using previously computed solutions as initial condition without the need to actually compute the tangent to the branch solution. It enables us nevertheless to reach high-energy solutions with very sharp gradients near the ball center; see Fig.~\ref{Fig_Gelfand}.
Of course many other types of cost functionals exploiting an {\it a priori} knowledge (radially decreasing solutions, etc.) could have been imagined and we do not claim for ``optimal choice.''  

In comparison, mesh-based methods in 3-D for computing  the bifurcation diagram directly from  \eqref{Eq_PDE}  would involve very important resources. We expect that 
the number of degrees of freedom (e.g.~grid size) needed here would be of order $10^6$.  In contrast, our results are obtained with NNs involving about 3000 degrees of freedom. Of course, in dimension $N\geq 4$, no mesh-based method is able to cope with such a problem. We address next this challenge within our neural network framework.

\begin{figure}[htbp]
\centering
\includegraphics[width=.9\textwidth,height=.6\textwidth]{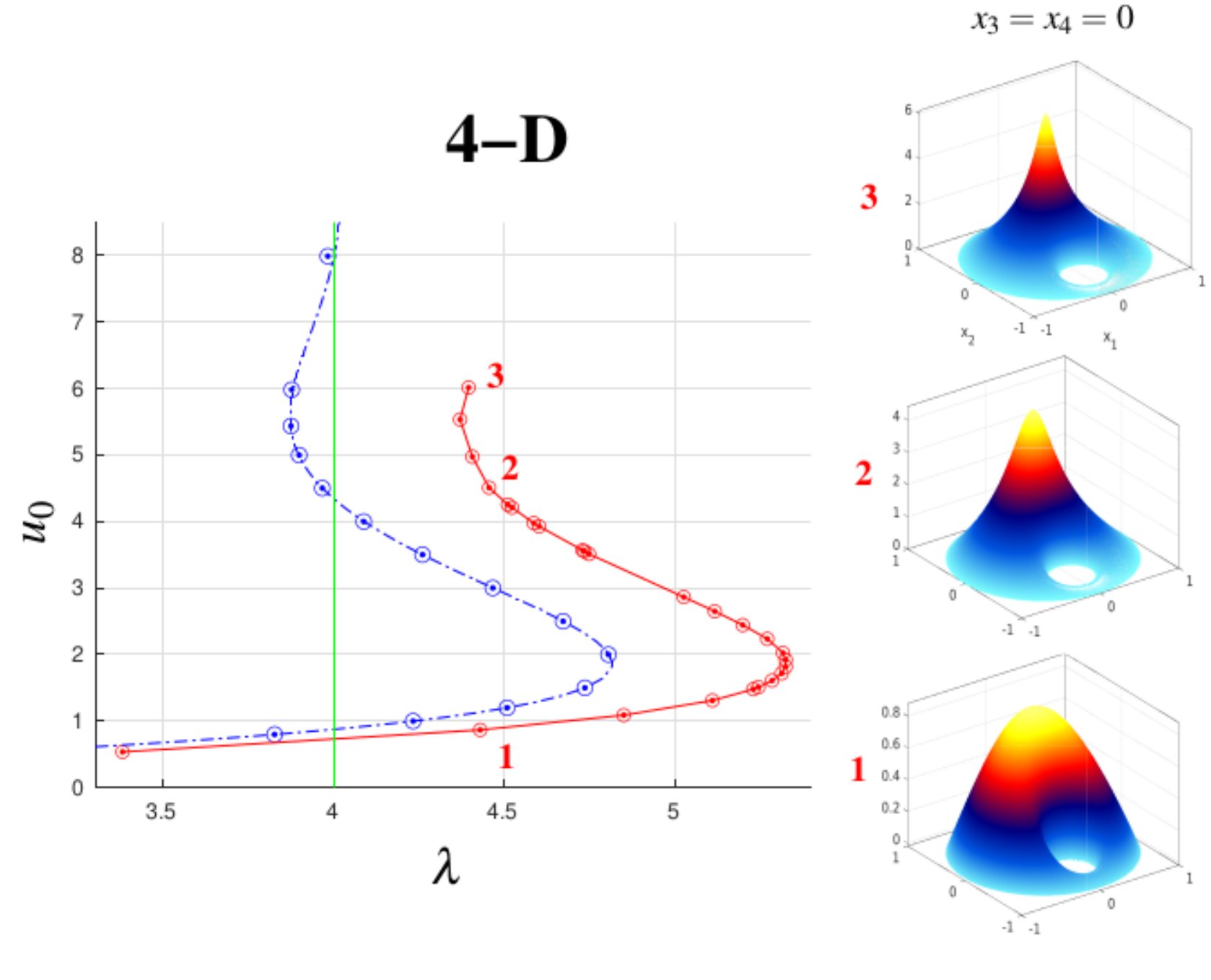}
	\caption{Bifurcation diagram in 4-D for the case $\Omega=B({\bf 0},1)\backslash B_r(\x_s)$, i.e.~the ball with a hole. This diagram is obtained by minimizing \eqref{Eq_initialC}. Here {\mk $\x_s\in B_{2/3}({\bf 0})$} and $r=1/4$.  The cross-sections in the plane $x_3=x_4=0$ of the solutions marked by 1,2, and 3 on this bifurcation diagram are shown in the left panels marked respectively by 1,2, and 3. As a comparison, the solutions obtained by minimizing \eqref{Eq_initialC} in the case of the (full) ball $\Omega=B({\bf 0},1)$ are represented by blue dots, while the blue curve is obtained by using a high resolution  pseudo-arclength code solving the 1-D problem \eqref{Eq_1D}.} 
\label{Fig_cyclope}
\end{figure}

\subsection{The case of domains with holes in dimension $N=4$}
In the previous case we demonstrated the ability to compute the bifurcation diagram by exploiting an a priori knowledge on the problem, using a transformation allowing for smoothing out the singular behavior near the ball's center. 
Here, we consider the unit ball in dimension 4 with cavities that break the radial symmetry and that thus prevent us to use such an an a priori knowledge. {\mk For these domain configurations} with $f(u)=e^u$, a few properties are known about the global shape of the bifurcation diagram, but not about its details.  The known features include the existence of a critical $0<\lambda^\ast<\infty $,  a branch of minimal solutions $u_\lambda^{\#}$ in $H_0^1(\Omega) \cap L^{\infty}(\Omega)$ such that $\lambda \mapsto u_\lambda^{\#}$ is increasing over $(0,\lambda^\ast)$, and that the full solution set $\{(\lambda,u_\lambda)\}$ forms an unbounded continuum in $[0,\lambda^\ast)\times C^2(\overline{\Omega})$; see \cite{rabinowitz1971some,amann1976fixed,lions1982existence,cazenave2006introduction} and \cite[Appendix A]{chekroun2018topological} for a self-contained expository of the latter point.  Very little is known however about the shape of the solutions that populate such a continuum and if the latter  has, as in the case of the ball in dimension four, many (possibly infinite) turning points.  

Our approach allows us to provide the first numerical hints in 4-D in the case $\Omega=B({\bf 0},1)\backslash B_r(\x_s)$, where $B_r(\x_s)$ is the closed ball centred at $\x_s$ in $B({\bf 0},1)$ of radius $r<1-\norm{\x_s}$. Our results show indeed that at least two turning points exist for this case and that, as in the case of the full ball, the solution becomes more and more singular as one ``climbs'' along the bifurcation diagram; see Fig.~\ref{Fig_cyclope}.

These results are obtained by minimizing \eqref{Eq_initialC} in which the Dirichlet boundary conditions are handled here again via a lift procedure which consists of using the ansatz  
$u(\x) = A \Theta(\x) v(\x)$ with $\Theta(\x) = \cos(\pi \|\x\|^2/2)$ for {\mk computing} the branch of minimal solutions $u_\lambda^\#$, and $\Theta(\x) = u_\lambda^\#(\x)$, after the first turning point has been crossed. The reason of changing of lift function after the first turning point is that it enables for encoding the sharp solutions' gradients that develop within $B({\bf 0},1)$ {\mk near the hole's boundary.}
By doing so, the {\mk NN-solver} is able to reach a mean-square error of $10^{-6}$ on the {\mk ``internal and external'' domain's boundaries} for the solutions shown in Fig.~\ref{Fig_cyclope}. 
A comparison with the case of the four-dimensional ball (without hole), shows that the bifurcation diagram in the case $\Omega=B({\bf 0},1)\backslash B_r(\x_s)$ shares a similar shape albeit {\mk with a first turning point} stretched to the right; compare blue and red curves in Fig.~\ref{Fig_cyclope}.

\begin{figure}[htbp]
\vspace{-.2cm}
\centering
\includegraphics[width=.7\textwidth,height=.4\textwidth]{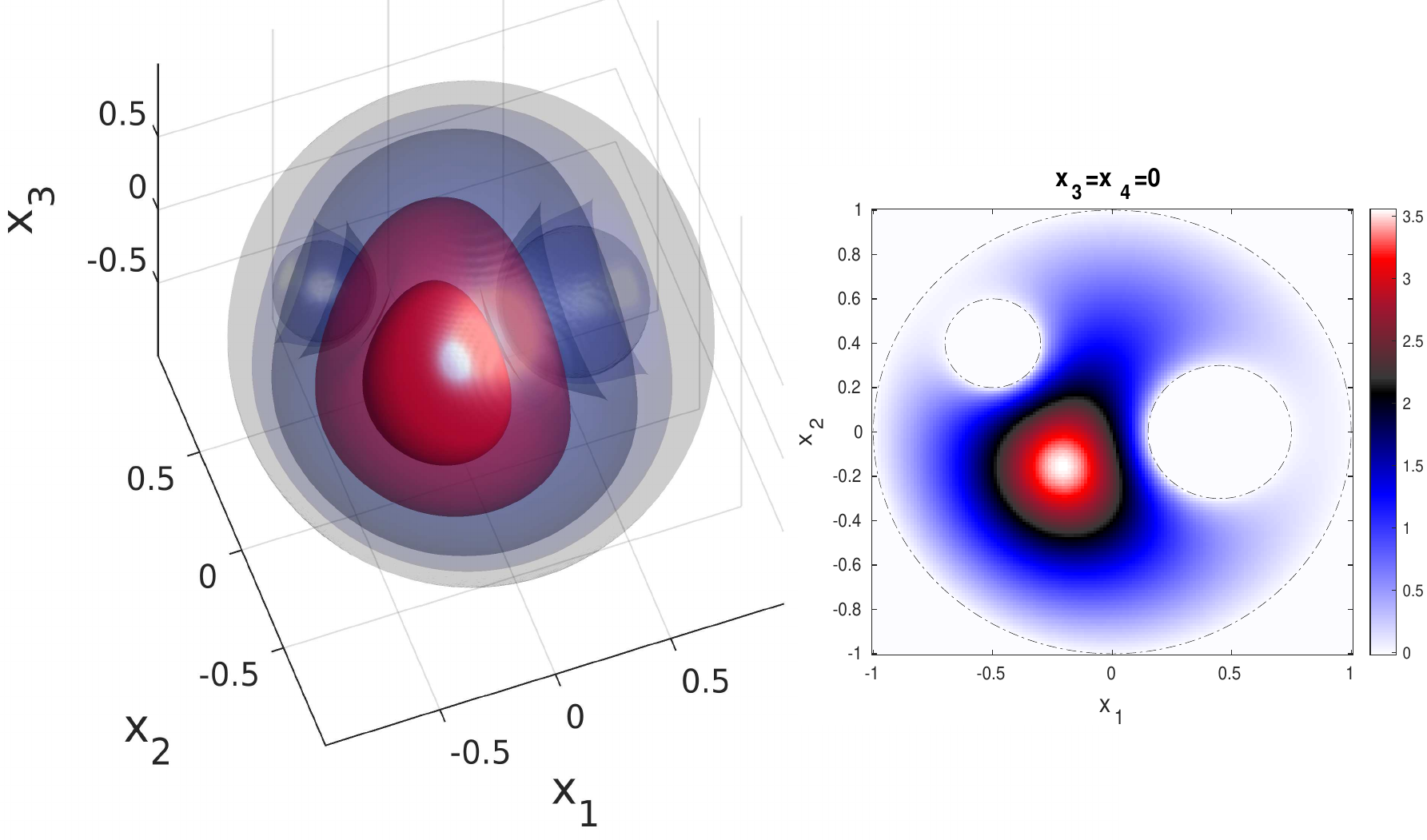}
	\caption{Three-dimensional (left panel) and a two-dimensional (right panel) sections of a solution to the Gelfand problem in 4-D over the unit ball with two cavities.}
\label{Fig_clown}
\end{figure}

Transformations inspired by \eqref{Eq_transfo} exploiting estimates about the location of the singularity could be used to reach out higher-energy solutions, but this requires more work. The approach is any way versatile enough to handle more complex geometries in dimensions higher than three. In that respect, Figure \ref{Fig_clown} shows {\mk 3-D and 2-D sections} of a solution to the Gelfand problem in 4-D over a domain with two holes of different size. It corresponds to a energy-level of type 2, i.e.~after the first turning point,  shown in  Fig.~\ref{Fig_cyclope} in the case of a single hole.

\section{Discussion and Outlook}
Thus, we provided a flexible machine learning framework  using simple FFNNs, geared toward high-dimensional eigenvalues problem of diffusion operators, linear and nonlinear, that are beyond reach with mesh-based methods. It opens up a wide range of applications and extensions for further research.

 For instance, in computer vision problems involving partial shape similarities, it is known that matching similar regions in 3D can be formulated as an alignment of $k$ eigenvalues of operators closely related to the Laplace-Beltrami operator (LBO)  \cite{rampini2019correspondence}. As the number $k$ is getting large ($k>100$) to favor a better shape discrimination,   high-precision discretization schemes of the LBO operator are however  required to avoid artifacts related to mesh tessellation; see \cite[Fig.~11]{rampini2019correspondence}. The NN-approach proposed here  allowing for the simultaneous computation of eigenvalues of such operators, could provide a natural mesh-free remedy to this problem. 
 
As mentioned above, the approach presented here is not limited to second-order differential operators. In that respect, the computation of bifurcation diagrams for Gelfand-type problems involving the $p$-Laplacian operator \cite{jacobsen2002liouville} in non spherical geometries and in higher dimensions could be addressed in a similar fashion. 

 Finally, in the vast topic of light scattering, vector Helmholtz equations with a small parameter  or discontinuous coefficients are known to play a prominent role (e.g.~polarization) \cite{gouesbet2011generalized}. Certain eigensolutions may exhibit very complicated shapes which require already in 2D intensive computations on a supercomputer \cite[Fig.~10]{gagnon2015lorenz} and are out of reach by the traditional series expansions used in the field \cite{gouesbet2011generalized,hergert2012mie}. We hope that addressing such Helmholtz problems within our NN-framework could provide an alternative approach for computing such eigenmodes with much less computational efforts, including in 3D.

In these problems or those considered in this work, the proper handling of eigenmodes' many possible sharp gradients over small regions is key to resolve. A natural idea for improving the performance of the proposed NN-approach consists of adaptively sampling more points in locations where the residual is large or use generative adversarial neural network to figure out where the NN is likely to be incorrect. We leave these important practical aspects for future investigations.


 \section*{Acknowledgments}
 
 This work has been supported by the Office of Naval Research (ONR) Multidisciplinary University Research Initiative (MURI) grant N00014-20-1-2023. This study was also supported by a Ben May Center grant for theoretical and/or computational research and by the Israeli Council for Higher Education (CHE) via the Weizmann Data Science Research Center, and by the European Research Council (ERC) under the European Union’s Horizon 2020 research and innovation programme (Grant Agreement No. 810370).

\appendix
\section{Material and Methods}
\small{
\subsection{Neural network parametrisation}
 Our neural network model is aimed at mapping the input vector $\x$ in $\R^{N}$ onto the vector $(\phi_1(\x),\cdots,\phi_k(\x))$, made of $k$ eigenmodes evaluated at $\x$. The sizes of the input and output features is respectively determined by the dimension $N$ of the ambient space and the number $k$ of eigenpairs $(\lambda,\phi)$ we target to approximate.

We denote our model output by $\boldsymbol{\Phi}(\x,\T)$, where $\T$ is the vector of network parameters including weights and biases. {\mk In the case of e.g.~the Kolmogorov operator,} these are found by solving the optimization problem  \eqref{Eq_multiopt} in the case $k>1$ and by minimizing the  cost functional \eqref{Eq_Cost} in the case $k=1$.   {\mk The derivatives are represented by finite differences on the NN parameterization.}

The neural network processes the input features using a number of layers, each of which combines basic operations such as affine transformations and element-wise nonlinearities.

There is flexibility in choosing the size of the hidden layers, which is also called their widths, while the number of layers is called the network depth. Altogether, the depth, width, and design of the layers are referred to as the network's architecture. In this article, we work with  standard feedforward neural networks (FFNNs).

Thus, our multivalued neural network parametrisation of $k$ eigenmodes to e.g.~\eqref{Eq_eigpb}, takes the form
\be\label{Eq_phi}
\Phi(\x;\T) = {\N}_{\rm out} \circ {\N}_L \circ \cdots \circ {\N}_1 \circ {\N}_{\rm in}(\x),
\ee
with ${\N}_k(\y ) = \sigma_k({\bf W}^{(k)} \y + {\bf b}^{(k)})$ and $\T = \{({\bf W}^{(k)},\boldsymbol{b}^{(k)}), \; 1\leq k\leq L \}$, $L$ denoting the number of layers. 
Here $\circ$ denotes the composition operator. 
The terms ${\bf W}^{(k)}$ and ${\bf b}^{(k)}$ are 
the aforementioned weight matrices and bias vectors and can have possible variable sizes. 
Figure \ref{Fig_totoNN} shows a schematic of such a standard FFNN architecture to approximate the solutions to e.g.~\eqref{Eq_eigpb}   via minimization of   \eqref{Eq_multiopt} in the case $k>1$.

 It is important to note that the use of the neural network renders the minimization problem \eqref{Eq_multiopt}  non-convex.
Due to this property, standard optimisation 
technics do not operate, even more so when the input dimension becomes large. 
This is the reason why stochastic gradient descent algorithms are so popular: there are
efficient in filtering out the fine-grained nonconvex structures, in a statistical sense. Here, a diagnostic
is taken only when the NNs have reached some statistical equilibrium. In essence, stochastic gradient descent
is simply the integration through many iterations of a system of the form $\dot \T = - 1/p \sum_{k=1}^p \nabla {\C}(\T,\xi_k)$ where
 $\xi_k$  are random points distributed according to a given probability measure,  ${\C}$ is the cost functional
to minimize, and $\T$ are the NN parameters.

\subsection{A quick guide for the practitioner}
 We provide here several remarks of practical importance. 
Unlike supervised, data-based approaches relying on large-dimensional input vectors, in our unsupervised, equation-based approach the dimension input is much smaller, typically in the range $O(1$--$100)$. Spatial derivatives in the course of the optimisation are approximated by second-order finite differences with a very small increment of size $\sqrt{\epsilon_{\textrm{machine}}}$, where $\epsilon_{\textrm{machine}}$ is the machine precision. The derivatives of  NN-parameters are computed with automatic differentiation. 
 The integrals are computed via simple Monte-Carlo empirical means:
$\int_{\M} f({\bf x})~\d{\bf x} \approx 1/N \sum_{k=1}^N f({\bf x}_k)$ where
${\bf x}_k$ are sampled uniformally inside the domain. Other strategies are possible including 
adaptive Monte-Carlo (VEGAS), \cite{vegas}, cubature and quasi-Monte-Carlo formulas \cite{qmc}.

 The optimization of the NN is prescribed by the choice of ad-hoc hyperparameters. The better this choice, the more efficient the optimization. There are three classes of hyperparameters which most often couple together in complicated and sometimes unexpected ways. The first one is
the structure of the neural network, which depth, which capacity, which activation functions? In our case,
since we deal with the most simple building blocks (FFNNs), this question is mostly related
to the depth and capacity of the network. The second class is related to the descent algorithm, in particular
the training rate and the batch size and more generally, properties of the descent algorithm itself.  The third
important class is composed of all the penalty parameters.
\begin{figure}[htbp]
\centering
\includegraphics[width=.6\textwidth,height=.35\textwidth]{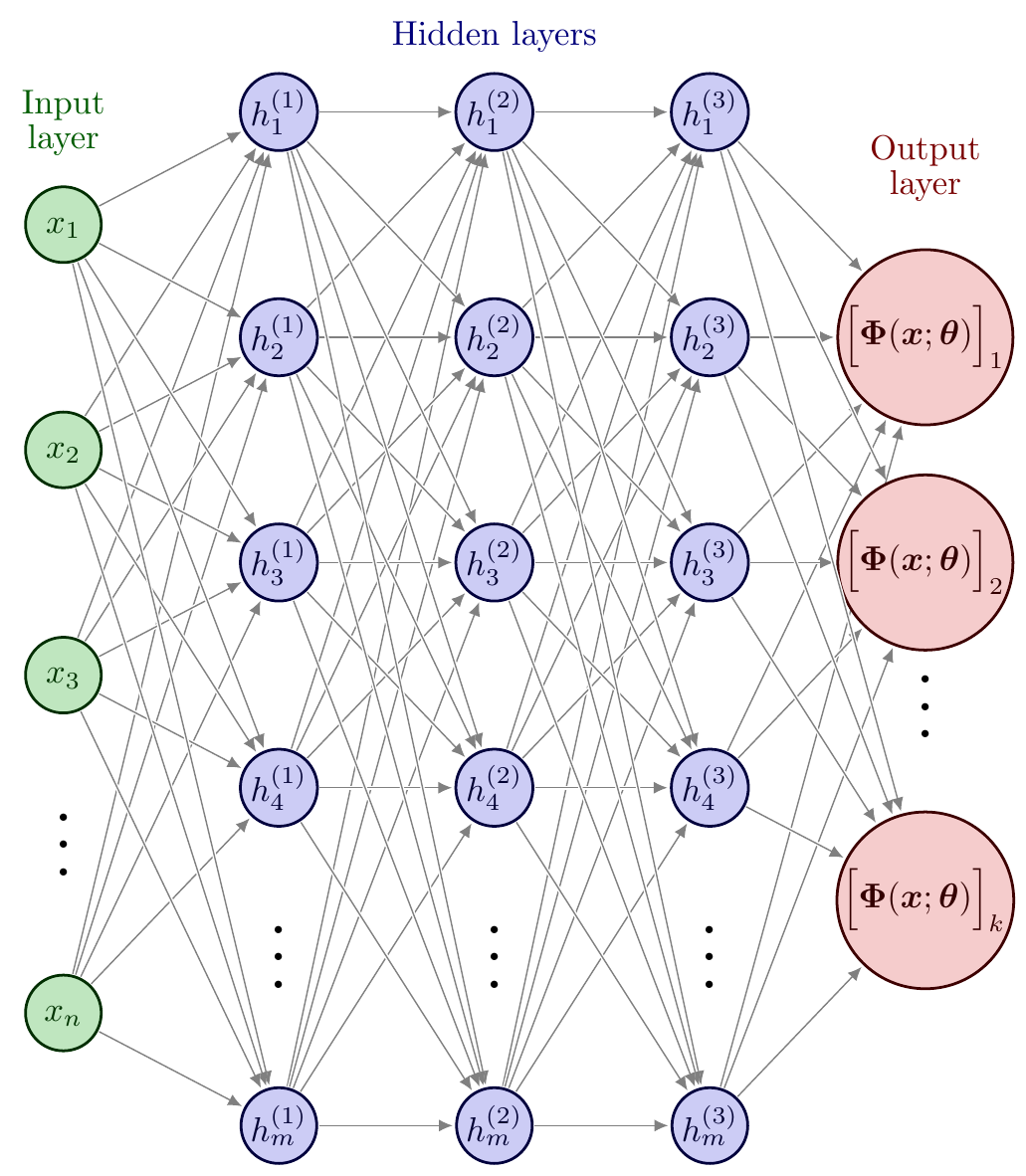}
	\caption{Schematic of the FFNN architecture used to minimize \eqref{Eq_multiopt}, i.e.~to deal with the simultaneous approximation of $k$ eigenmodes to \eqref{Eq_eigpb} in dimension $N$. Here the input is a point $\boldsymbol{x}=(x_1,\cdots,x_N)$ in $\mathbb{R}^N$ and three hidden layers are represented, each containing $m$ neurons.  The $k$-dimensional output $\boldsymbol{\Phi}(\x,\T)$ is aimed at approximating $k$ distinct eigenmodes evaluated at $\x$, i.e.~$(\phi_1(\x),\cdots,\phi_k(\x))$. Such an architecture is actually used in practice for the real and imaginary parts of the targeted eigenmodes.}
\label{Fig_totoNN}
\end{figure}

\begin{itemize}
\item[$\bullet$] {\bf Training rate}. In case it is too small, the descent may occur at an undesirable, too slow pace. On the other hand, a training rate that is too large drives the NN to display some spurious behaviours (e.g.~blow-up). One can this way identify relatively easily a reasonable range of training rates.
\item [$\bullet$] {\bf Batch size}. The rule of thumb is that the larger it is, the smaller is the variance. Dealing with large batches
has the drawback to be time consuming. Moreover, one can easily be trapped in unwanted regions of the landscape {\mk due to the small variance}.
However, in higher dimension, too small batches have the drawback to make the descent too slow.
\item[$\bullet$] {\bf Penalty coefficient}. The larger it is, the stiffer is the corresponding part of the cost functional which in turn imposes to use a smaller training rate. In some situations, it can be tricky to find some relevant ad-hoc range of values as particularly encountered to handle boundary conditions. This is one of the reasons why a lift procedure is often preferred at this stage; see Gelfand problem.  

\item[$\bullet$] {\bf Capacity}. In principle, the larger it is, the more accurate is the NN. In practice, however, we tend to favour
deep NNs with a ``narrow'' capacity. 
In our experiments, we observed that  using a capacity
larger by an order than the effective dimension of the problem is often helpful.
\end{itemize}

\subsection{The 10-D embedded stochastic Stepanov flow} 
The 10-D embedded stochastic Stepanov flow on the 10-dimensional torus is written as:
\bea\label{10dflow}
	\dot x_1 &= 
	 \alpha (1-\cos(x_1-x_6))+(1-\alpha)(1-\cos x_6) + \epsilon \d W_t^1\\
	 \dot x_2  &=  \sin x_7 \sin x_{10} + 2\cos 2x_4 - \cos(x_3-x_5+x_1)+ \epsilon \d W_t^2\\
\dot x_3 & =   \cos x_8 \sin(x_4+x_7)\hspace{-.5ex} - \hspace{-.5ex}3\cos 3x_5\hspace{-.5ex} +\hspace{-.5ex} 3\sin(x_1-2x_6)\hspace{-.5ex}+ \hspace{-.5ex}\epsilon \d W_t^3\\
\dot x_4 & =  \sin x_5 \sin(x_7-x_6) + 2\cos 2x_4 \hspace{-.5ex}+ \hspace{-.5ex}\cos x_{10} \cos x_1\hspace{-.5ex}+ \hspace{-.5ex}\epsilon \d W_t^4\\
\dot x_5 & =  -2 -2\cos x_4 \sin 2x_2 + \cos^2 x_5 + \cos 2x_8 + \epsilon \d W_t^5\\
\dot x_6 & =  \alpha(1-\cos(x_1-x_6))+ \epsilon \d W_t^6\\
\dot x_7 & =  3\cos x_7 \sin 2x_3 - \cos^2 x_2 + \cos(x_2-x_1)+ \epsilon \d W_t^7\\
	\dot x_8 & =  1 + \cos x_9 \sin 2x_{10} - \cos^2(x_8-x_3)+ \epsilon \d W_t^8 \\
\dot x_9 & =  \cos x_7 \sin 2x_6 - \cos^2 x_8 + 2\sin x_3+ \epsilon \d W_t^9\\
\hspace{-1ex} \dot x_{10} & =  \sin(2x_3\hspace{-.25ex}-\hspace{-.25ex}x_2)\hspace{-.5ex} - \hspace{-.5ex}\cos(x_{10}\hspace{-.25ex}-\hspace{-.25ex}x_6+\hspace{-.25ex}x_1)\hspace{-.5ex} + \hspace{-.5ex}\sin x_1 \hspace{-.5ex}\cos x_6\hspace{-.5ex}+ \hspace{-.5ex}\epsilon \d W_t^{\hspace{-.3ex}10}\\
\eea
where the $W_t^j$ are  mutually independent standard Brownian motions and $\alpha, \epsilon>0$.

\subsection{Gelfand problem in the radial case: Joseph \& Lundgren treatment}
We recall here the phase plane method of \cite{joseph1973quasilinear} to study the bifurcation diagram of the nonlinear eigenvalue problem \eqref{Eq_1D} allowing us to conclude easily to the existence of infinitely many turning points approaching a singular solution. 
First, consider the associated initial value problem (IVP) $u''=\frac{N-1}{r} u'+\lambda e^u$, where $u(0)=a$ and $u'(0)=0$,  with, $a$, a free parameter to be found such that $u(1)=0$.  Then, every solution to this IVP is obtained via the following Emden's transformation:
\bes
u(r)=w(t)-2t+a, \; r=\sqrt{\frac{2(N-2)}{\lambda e^a}}e^t,
\ees
in which $w$ solves $w'' + (N-2) w'+ 2(N-2) (e^w-1)=0$.
A phase plane analysis of this problem reveals that $(w,w')={\bf 0}$ is a stable focus for $3\leq N \leq 9$, with eigenvalues given by $2 \mu=-2-N\pm i\sqrt{(N-2)(10-N)}$. Note that $u(1)=0$ translates to $w(\tau)-2\tau +a=0$ 
which is equivalent to find $\tau$ such that $w(\tau)=\ln ( \lambda /2(N-2))$. Now since the orbit $\mathcal{O}(t)=(w(t),w'(t))$ is spiralling towards ${\bf 0}$,  we find at least $k$ solutions (for any $k$) for $\lambda$ close enough to $2(N-2)$, and infinitely many when $\lambda=2(N-2)$.

The orbit $\mathcal{O}(t)$ crosses the $w$-axis infinitely many times. One denotes by $\tau_k$ the crossing times for which $w'(\tau_k)=0$. Then $w(t)$ achieves either a local maximum or minimum at $t=\tau_k$, and  the $T_k=(\lambda_{\tau_k},u_{\lambda_{\tau_k}})$ corresponds to the turning points mentioned in the Main Text.

\subsection{Gelfand problem: Neural network configuration}
 In the case $\Omega=B(0,1)$ for $N=3$, the minimization of \eqref{Eq_rev_cost} following steps (G$_1$) and (G$_2$) has been operated by means of FFNNs with 12 hidden layers with 15 neurons for each layer giving rise to 2956 free parameters.   The activation functions
are {\mk swish} \cite{ramachandran2017searching} except for the output layer. The batch size has 512 points uniformly distributed on the 3-D unit ball (after the
blow-up transformation \eqref{Eq_transfo}.  The descent is executed using ADAM with a learning rate between $5 \cdot 10^{-4}$ and $10^{-4}$ depending on the value
of $A$: typically, the larger $A$, the smaller the training rate.  In step (G$_1$),  the penalisation coefficients are $\rho_g = 10$ 
and $\rho_{0} = 50$.  The parameter $A$ has been varied from $A=0.1$ to $A=14$ to obtain the results shown in Fig.~\ref{Fig_Gelfand}. The parameter $\epsilon$ in \eqref{Eq_transfo} is chosen to be $\epsilon=10$ on the branches below the second turning point, $\epsilon=10^{-2}$ on the branch right below the third one, and $\epsilon=10^{-3}$, after.

\subsection{Gelfand problem: Radially-scaled change of coordinates}
{\mk We provide here for the reader's convenience, the change-of-variable formulas used for transforming the Gelfand problem.  By introducing ${\bf X} = {\bf X(\x)}$, with ${\bf X(\x)}$ sufficiently smooth, we have trivially 
$$
\Delta_{\bf X} = \sum_{k=1}^N \Delta X_k \partial_k + \sum_{i,j =1}^N \left( \nabla X_i \cdot \nabla X_j \right) \partial_{ij},
$$
where $\Delta$ and $\nabla$ are taken with respect to $\x$. 

Assume that 
$$
{\bf X} (\x) = \frac{g(r)}{r} \x, \mbox{ with } r = \norm{\x},
$$ 
with $g$ some smooth invertible function of the real line. 
In this case, the inverse transformation is trivially given by 
$$
\x = \frac{g^{-1}(R)}{R} {\bf X}, \mbox{ with } R = \norm{{\bf X}}.
$$

Then, by introducing $R= g(r)$, $r = f(R)$  (i.e.~$f=g^{-1}$), we get after simplifications that  
\be
\nabla X_i \cdot \nabla X_j = \left( \frac{1}{R^2 (f'(R))^2} -\frac{1}{f^2(R)} 
\right) X_i X_j + \delta_{ij} \frac{R^2}{f^2(R)}.
\ee

Similarly, by expressing the 2nd-order derivatives of the variable $X_i = (g/r) x_i$, we arrive at
\bes
\Delta X_i = \frac{g''}{r} x_i + \frac{N-1}{r} \left( \frac{g}{r} \right)' x_i.
\ees

Thus, we have that: 
\bi
\item In coordinates $x_i$ with $()'$ denoting $\d/\d r$: 
$$
 \left\{
\begin{array}{llll}
\Delta X_i  & = & \Frac{g''}{r} x_i + \Frac{N-1}{r} \left(\Frac{g}{r} \right)' x_i \\\\
\nabla X_i \cdot \nabla X_j & = & \left( \Frac{g'^2}{r^2} - \Frac{g^2}{r^4} 
\right)  x_i x_j + \delta_{ij} \Frac{g^2}{r^2},
\end{array}
\right.
$$
with $r  =  g^{-1}(R)=f(R)$ and $x_i  = (f(R)/R)  X_i$. 

\item In coordinates $X_i$ with $()'$ denoting $\d/\d R$, we have: 
$$
\left\{
\begin{array}{llll}
\Delta X_i  & = & -\Frac{f''}{f'^3 R} X_i + \Frac{N-1}{R f'} \left(\Frac{R}{f} \right)' X_i \\\\
\nabla X_i \cdot \nabla X_j & = & \left( \Frac{1}{f'^2 R^2} - \Frac{1}{f^2} 
\right)  X_i X_j + \delta_{ij} \Frac{R^2}{f^2}.
\end{array}
\right.
$$
\ei 
}

}

\bibliographystyle{amsalpha}
\newcommand{\etalchar}[1]{$^{#1}$}
\providecommand{\bysame}{\leavevmode\hbox to3em{\hrulefill}\thinspace}
\providecommand{\MR}{\relax\ifhmode\unskip\space\fi MR }
\providecommand{\MRhref}[2]{%
  \href{http://www.ams.org/mathscinet-getitem?mr=#1}{#2}
}
\providecommand{\href}[2]{#2}

\end{document}